\documentclass[microtype]{agtart_a}
\pdfoutput=1
\usepackage{amscd}
\usepackage[all]{xy}

\title{Classifying spectra of saturated fusion systems}

\author{K\'ari Ragnarsson}
\givenname{K\'ari}
\surname{Ragnarsson}
\address{Department of Mathematical Sciences\\
University of Aberdeen\\\newline
Aberdeen AB24 3UE\\UK}
\email{kari@maths.abdn.ac.uk}



\received{25 March  2005}
\revised{19 January 2006}
\accepted{26 January 2006}
\proposed{}
\seconded{}
\publishedonline{26 February 2006}
\published{26 February 2006}

\volumenumber{6}
\issuenumber{}
\publicationyear{2006}
\papernumber{7}
\startpage{195}
\endpage{252}

\doi{}
\MR{}
\Zbl{}

\let\xysavmatrix\xymatrix
\def\xymatrix{\disablesubscriptcorrection\xysavmatrix}
\AtBeginDocument{\let\bar\wbar}
\allowdisplaybreaks
\newcommand{\Hom}[3]{\ensuremath{\operatorname{Hom}_{#3}\left(#1,#2\right)}}
\newcommand{\Mor}[3]{\ensuremath{\operatorname{Mor}_{#3}\left(#1,#2\right)}}
\newcommand{\Aut}[2]{\ensuremath{\operatorname{Aut}_{#2}\left(#1\right)}}
\newcommand{\Out}[2]{\ensuremath{\operatorname{Out}_{#2}\left(#1\right)}}

\newcommand{\Inj}[2]{\ensuremath{\operatorname{Inj}\left(#1,#2\right)}}

\DeclareMathAlphabet\EuR{U}{eur}{m}{n}
\SetMathAlphabet\EuR{bold}{U}{eur}{b}{n}
\newcommand{\curs}{\EuR}

\newcommand{\SpectrumCat}{\curs{Spectra}}
\newcommand{\GroupCat}{\curs{Gr}}

\newcommand{\Fp}{\ensuremath{\mathbb{F}_p}}
\newcommand{\Zp}{\ensuremath{\pComp{\mathbb{Z}}}}

\newcommand{\Coh}[1]{\ensuremath{H^{*}(#1)}}
\newcommand{\pComp}[1]{\ensuremath{{#1}^{\wedge}_p}}
\newcommand{\MappingClasses}[2]{\ensuremath{[#1,#2]}}

\newcommand{\F}{\ensuremath{\mathcal{F}}}
\newcommand{\HomF}[2]{\ensuremath{\operatorname{Hom}_{\F}\left(#1,#2\right)}}
\newcommand{\N}[4]{\ensuremath{ \{ #4 \in N_{#3}(#2) \mid #1 \circ c_{#4} \circ #1^{-1} \in Aut_{#3}(#1(#2)) \} }}  %% map, source, target, element
\newcommand{\CharIdem}{\ensuremath{\omega}}
\newcommand{\StableCharIdem}{\ensuremath{\wtilde{\omega}}}

\newcommand{\PresFunctor}[1]{\ensuremath{F_{#1}}}
\newcommand{\FusSysCat}{\curs{FS}}
\newcommand{\SatFusSysCat}{\curs{SFS}}

\newcommand{\ClSp}{\pComp{|\Link|}}
\newcommand{\Link}{\ensuremath{\mathcal{L}}}
\newcommand{\plfg}{\ensuremath{\left(S,\F,\Link\right)}}

\newcommand{\Stable}[1]{\Sigma^{\infty}{#1}}
\newcommand{\PtdStable}[1]{\Sigma^{\infty}_{+}{#1}}
\newcommand{\StableMaps}[2]{\{{#1},{#2}\}}
\newcommand{\PtdStableMaps}[2]{\{{#1}_+,{#2}_+\}}
\newcommand{\SphereSpectrum}{\ensuremath{\mathbb{S}^0}}

\newcommand{\Morita}{A^{+}}
\newcommand{\countS}{\ensuremath{\epsilon}}

\newcommand{\lsup}[2]{\ensuremath{{\vphantom{#1}}^{#2}\negmedspace{#1}}}
\newcommand{\DCF}[5] {\ensuremath{\sum_{x \in #1\backslash #5 /#4(#3)} \DCFterm{#1}{#2}{#3}{#4} }} % P, \varphi, Q, \psi, S
\newcommand{\DCFterm}[4] {\ensuremath{\left[{#4}^{-1}\left(#4\left(#3 \right) \cap #1^{x} \right),#2 \circ c_x \circ #4  \right]}}

\newcommand{\ClSpectrum}[1]{\ensuremath{\mathbb{B}{#1}}}
\newcommand{\ClSpecMaps}[2]{\ensuremath{[\ClSpectrum{#1},\ClSpectrum{#2}]}}
\newcommand{\StrClSpecFunctor}{\ensuremath{\Upsilon}}
\newcommand{\ClSpectraMap}[3]{\ensuremath{\ClSpectrum{#1_{#2}^{#3}}}}
\newcommand{\SCSTargetCat}{\ensuremath{\mathcal{C}}}

\newcommand{\Transfer}[1]{\ensuremath{Tr\left(#1\right)} }

\newcommand{\TransferFunctor}{\ensuremath{\Xi}}

\newcommand{\TransferTargetCat}{\ensuremath{\mathcal{D}}}
\newcommand{\InjSatFusSysCat}{\SatFusSysCat_{0}}

\newcommand{\Graph}[2]{\ensuremath{\Delta_{#1}^{#2}}}

\newcommand{\Fsub}[1]{\ensuremath{\mathop{\precsim}\limits_{#1}}}
\newcommand{\Fstrictsub}[1]{\ensuremath{\mathop{\precnsim}\limits_{#1}}}
\newcommand{\Fcon}[1]{\ensuremath{\mathop{\sim}\limits_{#1}}}
\newcommand{\Ffcon}[1]{\ensuremath{\mathop{\hbox{\small$\sim$}}\limits_{#1}}}

\newcommand{\noFsub}{\ensuremath{\precsim}}
\newcommand{\noFstrictsub}{\ensuremath{\precnsim}}
\newcommand{\noFcon}{\ensuremath{\sim}}
\newcommand{\noFfcon}{\ensuremath{\hbox{\mathsurround 0pt\small$\sim$}}}

\newcommand{\FsubModule}[4]{\ensuremath{M\left(\noFsub [#3,#4], #1, #2 \right)}}
\newcommand{\FconModule}[4]{\ensuremath{M\left(\noFcon [#3,#4], #1, #2 \right)}}
\newcommand{\FstrictsubModule}[4]{\ensuremath{M\left(\noFstrictsub [#3,#4], #1, #2 \right)}}

\newcommand{\noFsubModule}[2]{\ensuremath{M\left(\noFsub [#1,#2] \right)}}
\newcommand{\noFconModule}[2]{\ensuremath{M\left(\noFcon [#1,#2] \right)}}
\newcommand{\noFstrictsubModule}[2]{\ensuremath{M\left(\noFstrictsub [#1,#2] \right)}}

\newcommand{\StableFusion}[2]{\ensuremath{\pComp{A_{#1}\!\left(#2,#2\right)}}}

\def\HoCoPfeil{{\unitlength 1em\begin{picture}(0,.1)
\put(0,.1){\vector(1,0){3.75}}\end{picture}}}
\def\Hocolimname{{\unitlength.1em
\raisebox{-2.7\unitlength}{\begin{picture}(37,9.5)(0,0)
\put(0,2.7){$\operatorname{HoColim}$} \put(.1,-.1){\HoCoPfeil}
\end{picture}}}}
\def\Hocolim{\mathop{\Hocolimname}}

\newcommand {\hocolim}[1] {\Hocolim_{#1}}

\def\LinksPfeil{{\unitlength 1em\begin{picture}(0,.1)
\put(1.4,.1){\vector(-1,0){1.5}}
\end{picture}}}
\def\invlimname{{\unitlength.1em
\raisebox{-2.7\unitlength}{\begin{picture}(15.5,9.5)(0,0)
\put(0,2.7){$\operatorname{lim}$} \put(.05,-.1){\LinksPfeil}
\end{picture}}}}
\def\Invlim{\mathop{\invlimname}}

\newcommand {\invlimlim}[1] {\Invlim\limits_{#1}}

\makeatletter
\def\cnewtheorem#1[#2]#3{\newtheorem{#1}{#3}[section]
\expandafter\let\csname c@#1\endcsname\c@theorem}

\newtheorem{theorem}{Theorem}[section]
\cnewtheorem{proposition}[theorem]{Proposition}
\cnewtheorem{lemma}[theorem]{Lemma}
\cnewtheorem{corollary}[theorem]{Corollary}
\cnewtheorem{fact}[theorem]{Fact}
\newtheorem{maintheorem}{Theorem}

\theoremstyle{definition}
\cnewtheorem{definition}[theorem]{Definition}

\newtheorem*{question*}{Question}

\theoremstyle{definition}
\cnewtheorem{remark}[theorem]{Remark}
\cnewtheorem{example}[theorem]{Example}

\makeautorefname{maintheorem}{Theorem}

\makeatother  %  move after \newtheorem block

\begin{document}

\begin{asciiabstract}
The assignment of classifying spectra to saturated fusion systems
was suggested by Linckelmann and Webb and has been carried out by
Broto, Levi and Oliver. A more rigid (but equivalent) construction
of the classifying spectra is given in this paper. It is shown that
the assignment is functorial for fusion-preserving homomorphisms in
a way which extends the assignment of stable p-completed
classifying spaces to finite groups, and admits a transfer theory
analogous to that for finite groups. Furthermore the group of
homotopy classes of maps between classifying spectra is described,
and in particular it is shown that a fusion system can be
reconstructed from its classifying spectrum regarded as an object
under the stable classifying space of the underlying p-group.
\end{asciiabstract}

\begin{abstract}
The assignment of classifying spectra to saturated fusion systems
was suggested by Linckelmann and Webb and has been carried out by
Broto, Levi and Oliver. A more rigid (but equivalent) construction
of the classifying spectra is given in this paper. It is shown that
the assignment is functorial for fusion-preserving homomorphisms in
a way which extends the assignment of stable $p$--completed
classifying spaces to finite groups, and admits a transfer theory
analogous to that for finite groups. Furthermore the group of
homotopy classes of maps between classifying spectra is described,
and in particular it is shown that a fusion system can be
reconstructed from its classifying spectrum regarded as an object
under the stable classifying space of the underlying $p$--group.
\end{abstract}

\maketitle

\section*{Introduction}
Saturated fusion systems were introduced by Puig in \cite{Puig2,Puig} 
as a formalization of fusion systems of groups. To a finite
group $G$ with Sylow $p$--subgroup $S$ one associates a \emph{fusion
system $\F_S(G)$ over $S$.} This is the category whose objects are
the subgroups of $S$, and whose morphisms are the conjugations
induced by elements in $G$. Puig axiomatized this construction, thus
allowing abstract fusion systems without requiring the presence, or
indeed existence, of an ambient group $G$. He also identified
important properties enjoyed by those fusion systems that are
induced by groups. Puig called fusion systems with these properties
\emph{full Frobenius systems}. These definitions were later
simplified by Broto--Levi--Oliver, who introduced the term
\emph{saturated fusion systems} in \cite{BLO2} (see 
\fullref{def:SatFS} below). A further simplification has been obtained
by Kessar--Stancu in \cite{KS}.

A useful tool for the study of saturated fusion systems would be a
functor assigning a classifying space to each saturated fusion
system. Exactly what a classifying space means in this context is
made precise by the theory of $p$--local finite groups developed by
Broto--Levi--Oliver in \cite{BLO2}. They define a $p$--local finite
group as a triple $(S,\F,\Link)$, where $S$ is a finite $p$--group,
$\F$ is a saturated fusion system over $S$, and $\Link$ is a centric
linking system associated to $\F$, a category which contains just
enough information to construct a classifying space $\ClSp$ for
$\F$.

The motivating example for the definition of a $p$--local finite
group comes from finite groups. In \cite{BLO1}, Broto--Levi--Oliver
give an algebraic construction for a centric linking system
$\Link_S^c(G)$ associated to the fusion system $\F_S(G)$ of a finite
group $G$, and show that \mbox{$\pComp{|\Link_S^c(G)|} \simeq
\pComp{BG}$}.

Given the classifying space $\ClSp$, one can by \cite{BLO2}
reconstruct the fusion system via the following homotopy-theoretic
construction:
$$\HomF{P}{Q} = \{\varphi \in \Hom{P}{Q}{}  \mid \theta \circ B\iota_Q \circ B\varphi \simeq \theta \circ B\iota_P \},$$
where $\iota_P$ and $\iota_Q$ are the inclusions of the subgroups
$P$ and $Q$ in $S$, and $\theta$ is the natural ``inclusion''
\mbox{$BS \to \ClSp$}. This construction was first applied by
Martino--Priddy in \cite{MP} to show that if the $p$--completed
classifying spaces of two finite groups have the same homotopy
types, then their fusion systems are isomorphic.

The passage from saturated fusion systems to classifying spaces is
more problematic. In general it is not known whether a saturated
fusion system has an associated centric linking system, and if so,
whether it is unique. Broto--Levi--Oliver have developed an
obstruction theory to address these questions of existence and
uniqueness. Oliver has shown in \cite{Ol1,Ol2} that these
obstructions vanish for fusion systems of finite groups. Therefore
$\Link_S^c(G)$ is, up to equivalence, the unique centric linking
system associated to the fusion system $\F_S(G)$ of a finite group
$G$. Moreover, Oliver concludes that the $p$--completed classifying
spaces of two finite groups are homotopy equivalent if their fusion
systems over chosen Sylow subgroups are isomorphic via a
fusion-preserving isomorphism of these Sylow subgroups, thus proving
the Martino--Priddy conjecture \cite{MP}.

By Oliver's result the fusion system $\F_S(G)$ of a finite group
$G$ has a unique associated centric linking system. But even when
we restrict our attention to fusion systems coming from groups, we
do not have an expedient method to reconstruct the linking system
$\Link_S^c(G)$ (and consequently $\pComp{BG}$) from the fusion
data $\F_S(G)$. Nor do we know whether this assignment is
functorial, that is whether a morphism between fusion systems of
groups induces a map between their $p$--completed classifying
spaces.

A classifying space functor is not yet within our reach, but the
stable analogue presents a more tractable problem. When calculating
the cohomology of a $p$--local finite group $(S,\F,\Link)$ in
\cite{BLO2}, Broto--Levi--Oliver construct a characteristic biset for
$\F$. This is an $(S,S)$--biset $\Omega$ with properties, suggested
by Linckelmann--Webb, that guarantee that the induced stable selfmap
of $BS$ is an idempotent in cohomology with $\Fp$--coefficients (see
\fullref{prop:LWImportance}). Broto--Levi--Oliver noted that
the stable summand of $\Stable{BS}$ induced by a characteristic
biset $\Omega$ is independent of the particular choice of $\Omega$,
and agrees with the suspension spectrum of the classifying space
$\ClSp$. Furthermore, they observed that the construction of
$\Omega$ depends only on the saturated fusion system $\F$ and not on
the centric linking system $\Link$, and that therefore the induced
summand $\ClSpectrum{\F}$ can be considered as a \emph{classifying
spectrum} for the saturated fusion system $\F$.

In this paper we take their idea further. We give a different
formulation of the construction of a classifying spectrum of a
saturated fusion system $\F$, which allows us to retain more
information associated to $\F$. More precisely, we refine the
construction of the biset $\Omega$ in \cite{BLO2} to produce an
idempotent $\StableCharIdem$ in $\StableMaps{BS}{BS}$ with the
following stable idempotent analogues of the Linckelmann--Webb
properties:
\begin{enumerate}
  \item[(a)]
$\StableCharIdem$ is a $\Zp$--linear combination of homotopy
classes of maps of the form\break \mbox{$\Stable{B\varphi} \circ tr_P$},
where $P$ is a subgroup of $S$, \mbox{$\varphi \in \HomF{P}{S}$}
and $tr_P$ denotes the reduced transfer of the inclusion
\mbox{$P\leq S$}.
  \item[(b1)]
For each subgroup \mbox{$P \leq S$} and each \mbox{$\varphi \in
\HomF{P}{S}$}, the restrictions\break \mbox{$\StableCharIdem \circ \Stable
B\iota_{P}$} and \mbox{$\StableCharIdem\circ{\Stable{B\varphi}}$}
are homotopic as maps \mbox{$\Stable{BP} \to \Stable{BS}$}.
  \item[(b2)]
For each \mbox{$P \leq S$} and each \mbox{$\varphi \in
\HomF{P}{S}$}, the compositions \mbox{$tr_P \circ \StableCharIdem$}
and \mbox{$tr_\varphi \circ \StableCharIdem$}, where $tr_{\varphi}$
is the reduced transfer of the monomorphism
\mbox{$\varphi\negmedspace: P \to S$}, are homotopic as maps
\mbox{$\Stable{BS} \to \Stable{BP}$}.
  \item[(c)] $\wtilde{\countS}(\CharIdem) = 1$, where \mbox{$\wtilde{\countS} \negmedspace : \StableMaps{BS}{BS} \to \Zp$} is a
  morphism of modules derived from an augmentation of $\PtdStableMaps{BS}{BS}$  (see \fullref{lem:StableAugmentation} and \fullref{sec:StableCharIdempotent}).
\end{enumerate}
We show that $\StableCharIdem$ is the unique idempotent in
$\StableMaps{BS}{BS}$ with these properties and that Property (b1)
characterizes morphisms in the fusion system $\F$. Therefore we
write $\StableCharIdem_{\F}$ and refer to $\StableCharIdem_{\F}$ as
the \emph{stable characteristic idempotent} of $\F$.

The homotopy type of the stable summand of $\Stable{BS}$ induced by
$\StableCharIdem_{\F}$ agrees with the homotopy type of the
classifying spectrum $\ClSpectrum{\F}$ constructed by
Broto--Levi--Oliver, so this construction offers nothing new in
itself. It is the careful study of the characteristic idempotent
which allows us to exercise control over its mapping telescope
$\ClSpectrum{\F}$. We refer to the structure map
\mbox{$\sigma_{\F}\negmedspace: \Stable{BS} \to \ClSpectrum{\F}$} of
the mapping telescope as the \emph{structure map of $\F$}, and when
regarded as an object under $\Stable{BS}$, we refer to the pair
$(\sigma_{\F},\ClSpectrum{\F})$ as the \emph{structured classifying
spectrum of $\F$.} The structure map $\sigma$ admits a
\emph{transfer map} $t_{\F}$, which is, up to homotopy, the unique
map \mbox{$\ClSpectrum{\F} \to BS$} such that \mbox{$t_{\F} \circ
\sigma_{\F} \simeq \StableCharIdem_{\F}$} and \mbox{$\sigma_{\F}
\circ t_{\F} \simeq id_{\ClSpectrum{\F}}$}.

The reward for taking this point of view is the following result,
which further justifies the use of the term ``classifying
spectrum''. It appears in the text as \fullref{thm:sigmaToF}.
\begin{maintheorem} \label{mthm:sigmaToF}
If $\F$ is a saturated fusion system over a finite $p$--group $S$,
then $\F$ can be recovered from its structured classifying
spectrum $(\sigma_{\F},\ClSpectrum{\F})$ by the following
homotopy-theoretic construction:
$$\HomF{P}{Q} = \{\varphi \in \Hom{P}{Q}{}  \mid \sigma_{\F} \circ \Stable{B\iota_Q} \circ \Stable{B\varphi} \simeq \sigma_{\F} \circ \Stable{B\iota_P} \},$$
where $\iota_P$ and $\iota_Q$ are the inclusions of the subgroups
$P$ and $Q$ in $S$.
\end{maintheorem}
By Martino--Priddy \cite[Example 5.2]{MP3} the fusion system can not be recovered
from the homotopy type of the classifying spectrum alone; it must be
regarded as an object under $\Stable{BS}$. When applied to fusion
systems of groups, this theorem gives an alternative stable
classification of $p$--completed classifying spaces of finite groups,
which is in some sense finer than the one in \cite{MP3}. Combined
with the Martino--Priddy conjecture, this shows that the unstable
$p$--completed classifying space of a finite group is determined by
the stable $p$--completed classifying space, regarded as an object
under the stable classifying space of its $p$--Sylow subgroup. This
matter is taken up in Ragnarsson \cite{KR:NewStableClassification}.

The central result in this paper, which allows us to conduct the
necessary analysis of characteristic idempotents, is the
calculation of an explicit $\Zp$--basis for the submodule
$$\StableCharIdem_{\F_2} \circ \StableMaps{BS_1}{BS_2} \circ \StableCharIdem_{\F_1} \subset \StableMaps{BS_1}{BS_2},$$
for saturated fusion systems $\F_1$ and $\F_2$ over finite
$p$--groups $S_1$ and $S_2$, respectively. This module is naturally
isomorphic to the group of stable maps between the classifying
spectra of the fusion systems involved and so we get the following
theorem, a more concise version of which appears later as 
\fullref{thm:ClSpectraMapBasis}, as an immediate consequence.

\begin{maintheorem}\label{mthm:ClSpectraMapBasis}
Let $\F_1$ and $\F_2$ be saturated fusion systems over finite
\mbox{$p$--groups} $S_1$ and $S_2$, respectively. Then the group of
homotopy classes of stable maps from $\ClSpectrum{\F_1}$ to
$\ClSpectrum{\F_2}$ is a free $\Zp$--module with one basis element
\mbox{$\sigma_{\F_2} \circ (\Stable{B\psi} \circ tr_P) \circ
t_{\F_1}$} for every conjugacy class of pairs $(P,\psi)$
consisting of a subgroup \mbox{$P \leq S_1$} and a nontrivial
homomorphism \mbox{$\psi\negmedspace: P \to S_2$}. Conjugacy here
means that $\F_1$--conjugacy is taken in the source and
$\F_2$--conjugacy is taken in the target.
\end{maintheorem}

When $\F_1$ and $\F_2$ are fusion systems of groups, this theorem
can be applied to give a new variant of the Segal conjecture
describing the group of homotopy classes of stable maps between
$p$--completed classifying spaces of finite groups. This discussion
is taken up in Ragnarsson \cite{KR:NewSegal}.

If $\F_1$ and $\F_2$ are saturated fusion systems over finite
$p$--groups $S_1$ and $S_2$, respectively, the obstruction to
restricting a homomorphism \mbox{$\gamma\negmedspace: S_1 \to
S_2$} to a map between classifying spectra respecting their
structure maps is the compatibility of $\gamma$ with the stable
characteristic idempotents. This compatibility is achieved when
$\gamma$ is a\break \emph{$(\F_1,\F_2)$--fusion-preserving homomorphism},
which means that $\gamma$ induces a functor
\mbox{$\PresFunctor{\gamma}\negmedspace: \F_1 \to \F_2$} such that
\mbox{$\PresFunctor{\gamma}(P) = \gamma(P)$} for all \mbox{$P \leq
S_1$} and \mbox{$\gamma\vert_{Q} \circ \varphi =
\PresFunctor{\gamma}(\varphi) \circ \gamma\vert_P$} for all
\mbox{$\varphi\in\Hom{P}{Q}{\F_1}$}. Letting $\SatFusSysCat$
denote the category whose objects are the saturated fusion systems
and whose morphisms are fusion-preserving homomorph\-isms, we get
the following result which follows from 
\fullref{thm:StrClSpecFunctor} in the text.

\begin{maintheorem} \label{mthm:ClSpecFunctor}
There is a classifying spectrum functor
$$\ClSpectrum\negmedspace: \SatFusSysCat \longrightarrow \SpectrumCat$$
acting on objects by sending a saturated fusion system to its
classifying spectrum and on morphisms by sending a
$(\F_1,\F_2)$--fusion-preserving morphism $\gamma$ to the map
$$\ClSpectraMap{\gamma}{\F_1}{\F_2} := \sigma_{\F_2} \circ \Stable{B\gamma} \circ t_{\F_1} \negmedspace: \ClSpectrum{\F_1} \longrightarrow \ClSpectrum{\F_2},$$
which satisfies
$$\ClSpectraMap{\gamma}{\F_1}{\F_2} \circ \sigma_{\F_1} \simeq \sigma_{\F_2} \circ \Stable{B\gamma}. $$
\end{maintheorem}

It is an important property of this functor that when
\mbox{$\gamma\negmedspace: S_1 \to S_2$} is the restriction of a
homomorphism \mbox{$\overline{\gamma}\negmedspace: G_1 \to G_2$}
to Sylow subgroups, the map \mbox{$\ClSpectrum{\F_{S_1}(G_1)} \to
\ClSpectrum{\F_{S_2}(G_2)}$} induced by $\gamma$ is equivalent to
the map \mbox{$\pComp{\Stable{BG_1}} \to \pComp{\Stable{BG_2}}$}
induced by $\overline{\gamma}$, as maps of objects under the
stable classifying spaces of their Sylow subgroups. This is proved
in \fullref{sec:Comparison}.

A monomorphism \mbox{$\gamma \negmedspace : S_1 \to S_2$} admits a
transfer map \mbox{$tr_{\gamma}\negmedspace: \Stable{BS_2} \to
\Stable{BS_1}$}, which restricts to a map of classifying spectra
that preserves transfer maps when $\gamma$ is fusion-preserving.
Collecting \fullref{thm:TransferFunctor} and Propositions
\ref{prop:ActsByMultiplication} and \ref{prop:TransfersFrobenius},
we get the following result.

\begin{maintheorem} \label{mthm:Transfers}
There is an assignment of a transfer map
$$\Transfer{\gamma_{\F_1}^{\F_2}} := \sigma_{\F_1} \circ tr_{\gamma} \circ t_{\F_2}  \negmedspace:  \ClSpectrum{\F_2} \to \ClSpectrum{\F_1},$$
to every $(\F_1,\F_2)$--fusion-preserving monomorphism
\mbox{$\gamma\negmedspace: S_1 \to S_2 $}. The assignment has the
following properties:
\begin{enumerate}
  \item[(i)]
$t_{\F_1} \circ \Transfer{\gamma_{\F_1}^{\F_2}} \simeq tr_{\gamma}
\circ t_{\F_2}.$
  \item[(ii)]
$\Transfer{\gamma_{\F_1}^{\F_2}} \circ
\Transfer{\gamma_{\F_2}^{\F_3}} \simeq
\Transfer{\gamma_{\F_1}^{\F_3}}$.
  \item[(iii)]
The composition $\ClSpectrum{\gamma_{\F_1}^{\F_2}} \circ~
\Transfer{\gamma_{\F_1}^{\F_2}}$ acts on
$\Coh{\ClSpectrum{\F_2;\Fp}}$ as multiplication by $|S_2|/|S_1|$.
  \item[(iv)]
The transfer $\Transfer{\gamma_{\F_1}^{\F_2}}$ satisfies the
Frobenius reciprocity relation
$$\Transfer{\gamma_{\F_1}^{\F_2}}^*\left(\ClSpectrum{\gamma_{\F_1}^{\F_2}}^*\left(x\right) \cdot y\right) = x \cdot \Transfer{\gamma_{\F_1}^{\F_2}}^*\left(y\right)$$
for \mbox{$x \in \Coh{\ClSpectrum{\F_2};\Fp}$} and \mbox{$y \in
\Coh{\ClSpectrum{\F_1};\Fp}$}.
\end{enumerate}
\end{maintheorem}

The motivation for the work in this paper comes from a question of
Miller regarding an alternative formulation of $p$--local finite
groups in terms of homotopy subgroup inclusions satisfying a certain
transfer property. The author has obtained a partial answer to this
question, but on the way to doing so, has discovered results about
classifying spectra of saturated fusion systems which are most
likely of interest to a wider audience than the original question,
and are therefore presented separately in this paper. These results
do not depend on centric linking systems, and to emphasize this we
mostly avoid mentioning centric linking systems in this paper.
Implications for $p$--local finite groups will be discussed in a
subsequent paper \cite{KR:Transfers-plfgs}, where Miller's question
will be addressed.

\medskip
\textbf{Notational conventions}\qua Throughout this paper, $p$ is a
fixed prime. Cohomology is always taken with $\Fp$--coefficients. For
a space $X$ we let $X_+$ be the pointed space obtained by adding a
disjoint basepoint to $X$, and we let $\pComp{X}$ denote the
Bousfield--Kan $p$--completion \cite{Monster}.

The category of finite groups and homomorphisms is denoted by
$\GroupCat$. For an element $g$ of a group $G$, we let $c_g$ denote
the conjugation \mbox{$x \mapsto gxg^{-1}$}. When $H$ is a subgroup
of $G$ we write $\lsup{H}{g}$ for the conjugate $c_g(H)$ and $H^g$
for the inverse conjugate \mbox{$c_g^{-1}(H) = g^{-1}Hg$}. For
subgroups $H$ and $K$ of $G$ we let $N_G(H,K)$ denote the
transporter
$$N_G(H,K) := \{g \in G \mid \lsup{H}{h} \leq K \},$$
and write
\begin{align*}
\Hom{H}{K}{G} &:= \{ c_g : H\to K \mid g \in N_G(H,K) \} \\
              &\phantom{:}= N_G(H,K)/C_G(H)
\end{align*}
for the set of homomorphisms from $H$ to $K$ induced by
conjugation in $G$.

The inclusion of a subgroup $H$ into a supergroup is denoted by
$\iota_H$, specifying the supergroup when there is danger of
confusion. For the convenience of the reader we use the letters
$S$, $P$ and $Q$ to refer to finite $p$--groups, while $G$ and $H$
refer to general finite groups. Moreover we use $\varphi$ to
denote homomorphisms belonging to fusion systems, while $\psi$ and
$\rho$ denote general homomorphisms.

All stable homotopy takes place in the homotopy category of spectra
which we denote by $\SpectrumCat$. A discussion of the stable
homotopy category can be found for example in \cite{Ad:Stable}. We
will use the shorthand notation
$$\PtdStable{X} := \Stable{(X_+)} $$
for the suspension spectrum of $X_+$. Since we often have cause to
work with stable $p$--completed classifying spaces, we adopt the
shorthand notation
$$\ClSpectrum{(-)} :=\pComp{\Stable{B(-)}},$$
regarded as functors
$$\GroupCat \longrightarrow \SpectrumCat.$$
As is usual, for spaces $X$ and $Y$ we let $\StableMaps{X}{Y}$
denote the group of homotopy classes of stable maps
\mbox{$\Stable{X} \to \Stable{Y}$}, and for spectra $E$ and $F$ we
let $\MappingClasses{E}{F}$ denote the group of homotopy classes of
(degree 0) maps \mbox{$E \to F$}. All homotopies are unpointed.

\medskip
\textbf{Overview}\qua In the first section we recall the definition of
saturated fusion systems. The second section treats Burnside modules
and the Segal conjecture relating them to stable maps between
classifying spaces of groups. In addition we develop some tools and
notation we will use throughout the paper. In 
\fullref{sec:Subconjugacy} we introduce the notion of fusion
subconjugacy. For fusion systems $\F_1$ and $\F_2$ over finite
$p$--groups $S_1$ and $S_2$, this gives a useful fusion-invariant
filtration of the Burnside module $A(S_1,S_2)$. In 
\fullref{sec:Idempotent} we assign a characteristic idempotent
$\CharIdem_{\F}$ in the $p$--completed double Burnside ring
$\pComp{A(S,S)}$ to a saturated fusion system $\F$ over $S$. In
\fullref{sec:Inherent} we perform a careful analysis of the
inherent properties of this idempotent, and in 
\fullref{sec:StableCharIdempotent} we interpret these results for the
stable idempotent $\StableCharIdem_{\F}$ of $\ClSpectrum{S}$ induced
by $\CharIdem_{\F}$. In \fullref{sec:ClSpectrum} we define the
classifying spectrum of $\F$ as the summand $\ClSpectrum{\F}$ of
$\Stable{BS}$ given by $\StableCharIdem_{\F}$, and prove that this
assignment is functorial. In \fullref{sec:Transfer} we develop
the theory of transfers for classifying spectra, and in 
\fullref{sec:Cohomology} we look at the behaviour of classifying spectra
and their transfers in cohomology. We conclude this paper in 
\fullref{sec:Comparison} by showing that the theory of classifying
spectra of saturated fusion systems developed here agrees with
existing theories of stable classifying spaces of saturated fusion
systems.

\medskip
\textbf{Acknowledgements}\qua I would like to thank my thesis advisor
Haynes Miller for suggesting the problem out of which this work grew
and for his enthusiasm and helpful advice during its progress. I also
thank Bob Oliver for his emailed suggestions for the proof of
convergence in \fullref{sec:Idempotent} and Ran Levi for many
lively and encouraging discussions on this subject.  Finally I thank
the referee on behalf of both the reader and myself for a very
accurate and helpful report which has improved the exposition and
clarity of the paper. The method of proof of the central result in
this paper is, at least subconsciously, inspired by Nishida's work in
\cite{Ni} and a preliminary version thereof.  

The author was supported by EPSRC grant GR/S94667/01 during part of
this work.

\section{Saturated fusion systems}\label{sec:Def&Term}
In this section we recall the definition of a saturated fusion
system. We begin by presenting the motivating example.
\begin{definition}
Let $G$ be a finite group with Sylow $p$--subgroup $S$. The
\emph{fusion system of $G$ over $S$} is the category $\F_S(G)$
whose objects are the subgroups of $S$, and whose morphisms are
the homomorphisms induced by conjugation in $G$:
\[
  \Hom{P}{Q}{\F_S(G)} = \Hom{P}{Q}{G}.
\]
\end{definition}
Puig \cite{Puig2,Puig} axiomatized this construction as follows.
\begin{definition} \label{def:FS}
A {\it fusion system} $\F$ over a finite $p$--group $S$ is a
category, whose objects are the subgroups of $S$, and whose
morphism sets $\HomF{P}{Q}$ satisfy the following conditions:
\begin{enumerate}
  \item[(a)] $\Hom{P}{Q}{S} \subseteq \HomF{P}{Q} \subseteq
\Inj{P}{Q}$ for all $P,Q \leq S$.
  \item[(b)] Every morphism in
$\F$ factors as an isomorphism in $\F$ followed by an inclusion.
\end{enumerate}
\end{definition}
From the definition it is clear that every fusion system over $S$
contains the fusion system $\F_S(S)$ of $S$. We denote this fusion
system by $\F_S$ for short.

Fusion systems at this level of generality are not particularly
useful or interesting, so we restrict to a certain subclass of
fusion systems introduced by Puig in \cite{Puig2}. Puig identified
important properties enjoyed by fusion systems of groups, and called
fusion systems with these properties \emph{full Frobenius systems}.
His definitions were later simplified by Broto--Levi--Oliver in
\cite{BLO2}, where they suggested the name \emph{saturated fusion
systems}. A further simplification has been obtained by
Kessar--Stancu in \cite{KS}.

We present the Broto--Levi--Oliver version below, but before stating
the definition, we need to introduce some additional terminology. We
say that two subgroups $P,P' \leq S$ are {\it $\F$--conjugate} if
they are isomorphic in $\F$. A subgroup $P \leq S$ is {\it fully
centralized in \F} if \mbox{$\vert C_S(P)\vert \geq \vert
C_S(P')\vert$} for every \mbox{$P' \leq S$} that is $\F$--conjugate
to $P$. Similarly $P$ is {\it fully normalized in \F} if
\mbox{$\vert N_S(P)\vert \geq \vert N_S(P')\vert$} for every
\mbox{$P' \leq S$} that is $\F$--conjugate to $P$.

\begin{definition} \label{def:SatFS}
A fusion system $\F$ over a $p$--group $S$ is {\it saturated} if
the following two conditions hold:
\begin{enumerate}
  \item[(I)] If $P \leq S$ is fully normalized in $\F$, then $P$ is
also fully centralized in $\F$, and $p$ does not divide the index of
$\Aut{P}{S}$ in $\Aut{P}{\F}$.
%\mbox{$\Aut{P}{S} \in Syl_p(\Aut{P}{\F})$.}
  \item[(II)] If $P \leq S$ and \mbox{$\varphi
\in \HomF{P}{S}$} are such that $\varphi (P)$ is fully
centralized, then $\varphi$ extends to \mbox{$\bar \varphi \in
\HomF{N_{\varphi}}{S}$}, where
$$ N_{\varphi} = \N{\varphi}{P}{S}{g}.$$
\end{enumerate}
\end{definition}

This definition is rather technical, and as the conditions in the
definition are not used explicitly in this paper, it suffices for
the reader to keep in mind that Condition I is a ``prime to $p$'' or
``Sylow'' property, analogous to the fact that the index of a Sylow
subgroup in a finite group is not divisible by $p$. Condition II is
a ``maximal extension property'' which (in a non-precise sense and
when combined with Condition I) can be thought of as an axiomatic
replacement of Sylow's Second and Third Theorems.

The role of saturated fusion systems in the theory of classifying
spectra developed in this paper is as follows. In 
\fullref{sec:Idempotent} we construct a characteristic idempotent
$\CharIdem$ for a fusion system $\F$ with a characteristic biset
$\Omega$. These objects are defined precisely in 
\fullref{sec:Idempotent}, and for now it suffices to say that
characteristic bisets are finite $(S,S)$--bisets with properties
stipulated by Linckelmann--Webb. The classifying spectrum of $\F$
is then constructed using $\CharIdem$ in 
\fullref{sec:ClSpectrum}. A construction of characteristic bisets for
saturated fusion systems is given by Broto--Levi--Oliver in
\cite{BLO2}. This allows us to develop the theory of classifying
spectra of saturated fusion system. But existence of a classifying
spectrum for a fusion system $\F$ depends only on the existence of
a characteristic biset for $\F$, and the properties of classifying
spectra follow from the Linckelmann--Webb properties without using
the saturation axioms. The theory therefore extends automatically
to all fusion systems that have characteristic bisets. It is an
interesting question whether the existence of a characteristic
biset for a fusion system $\F$ implies saturation of $\F$. The
author believes this is true, which is why the results in this
paper are only presented for saturated fusion systems.

\section{Burnside modules and the Segal conjecture}\label{sec:Burnside}
In this section we give a brief discourse about how stable maps
between classifying spaces of finite groups $G_1$ and $G_2$ are
related to $(G_1,G_2)$--bisets.

For finite groups $G_1$ and $G_2$, let $\Morita(G_1,G_2)$ be the
set of isomorphism classes of finite sets with a right
$G_1$--action and a free left $G_2$--action. The disjoint union
operation makes $\Morita(G_1,G_2)$ into a commutative monoid. We
denote the Grothendieck group completion by $A(G_1,G_2)$ and refer
to it as \emph{the Burnside module of $G_1$ and $G_2$}. The reader
should beware that this is not standard terminology. The group
structure of $A(G_1,G_2)$ is easy to describe. It is a free
abelian group with one generator corresponding to each transitive
$(G_1,G_2)$--biset. We proceed to describe and parametrize these
basis elements below.
\begin{definition} \label{def:Pair}
Let $G_1$ and $G_2$ be finite groups. A $(G_1,G_2)$--\emph{pair} is
a pair $(H,\psi)$ consisting of a subgroup \mbox{$H \leq G_1$} and
a homomorphism
$$\psi\negmedspace: H \to G_2.$$
We say that two $(G_1,G_2)$--pairs $(H_1,\psi_1)$ and
$(H_2,\psi_2)$ are $(G_1,G_2)$--\emph{conjugate} if there exist
elements \mbox{$g \in G_1$} and \mbox{$h \in G_2$} such that
\mbox{$c_g(H_1) = H_2$} and the following diagram commutes
\[
\begin{CD}
{H_1} @ > {\psi_1} >> {G_2} \\
@ V \cong V c_g V @ VV c_h V \\
{H_2} @> {\psi_2} >> {G_2.} \\
\end{CD}
\]
\end{definition}
\begin{remark}
Define the \emph{graph} of a $(G_1,G_2)$--pair $(H,\psi)$ by
$$\Graph{H}{\psi} := \{(h,\psi(h) \mid h \in H \} \leq G_1 \times G_2.$$
It is easy to check that $(G_1,G_2)$--pairs $(H_1,\psi_1)$ and
$(H_2,\psi_2)$ are $(G_1,G_2)$--\emph{conjugate} if and only if their
graphs are conjugate in \mbox{$G_1 \times G_2$}.
\end{remark}

We denote the $(G_1,G_2)$--conjugacy class of a $(G_1,G_2)$--pair
$(H,\psi)$ by $[H,\psi]_{G_1}^{G_2}$ or, when there is no danger
of confusion, just $[H,\psi]$. With a slight abuse of notation we
will also let $[H,\psi]_{G_1}^{G_2}$ (or $[H,\psi]$) denote the
basis element of $A(G_1,G_2)$ corresponding to the conjugacy class
of the $(G_1,G_2)$--pair $(H,\psi)$. Thus $[H,\psi]$ represents the
isomorphism class of the $(G_1,G_2)$--biset
$$G_2 \times_{(H,\psi)} G_1 := (G_2\times G_1)/\sim,$$
with the obvious right $G_1$--action and left $G_2$--action, where
$$(x,gy)\sim (x\psi(g),y)$$
for \mbox{$x \in G_2~,y \in G_1$} and \mbox{$g \in H$}.

Given three finite groups $G_1,G_2,$ and $G_3,$ we get a morphism
of monoids
$$-\circ- : \Morita(G_2,G_3) \times \Morita(G_1,G_2) \to \Morita(G_1,G_3)$$
$$(\Omega,\Lambda) \mapsto \Omega \circ \Lambda := \Omega \times_{G_2} \Lambda,
\leqno{\hbox{by}}$$
which extends to a bilinear map
\begin{equation} \label{eq:BurnsidePairing}
A(G_2,G_3) \times A(G_1,G_2) \to A(G_1,G_3).
\end{equation}
This pairing can be described in terms of the basis elements using
the double coset formula.
\begin{equation} \label{eq:DCF}
 [K,\rho]^{G_3}_{G_2} \circ [H, \psi]^{G_2}_{G_1} =
 \DCF{K}{\rho}{H}{\psi}{G_2}_{G_1}^{G_3}.
\end{equation}
We pay special attention to the simple case where \mbox{$K =
G_2$}, so $\rho$ and $\psi$ are composable morphisms. In this case
the double coset formula simplifies to
\begin{equation} \label{eq:DCFsimple}
 [G_2,\rho]^{G_3}_{G_2} \circ [H, \psi]^{G_2}_{G_1} = [H,\rho \circ
 \psi]^{G_3}_{G_1}.
\end{equation}
For a finite group $G$ the pairing of \ref{eq:BurnsidePairing}
makes $A(G,G)$ into a ring which we call the \emph{double Burnside
ring of $G$}. This should not be confused with the \emph{Burnside
ring} $A(G)$ \cite{Dieck}. The latter is the Grothendieck group
completion of the monoid of isomorphism classes of finite left
$G$--sets. As a $\Z$--module, $A(G)$ is free $\Z$--module with one
generator $[G/H]$ for each conjugacy class of subgroups \mbox{$H
\leq G$}. As a ring, the multiplicative structure on $A(G)$ is
induced by Cartesian product and linear extension.

Bisets relate to stable maps via the Becker--Gottlieb transfer
\cite{BG}. We recall some basic properties of transfers here, and
refer the reader to \cite{Ad:InfLoopSp} for a more thorough
discussion. Given a finite covering \mbox{$f \negmedspace : X \to
Y$}, where $Y$ is connected, Becker--Gottlieb constructed a stable
map \mbox{$tr_f \negmedspace : \PtdStable{Y} \to \PtdStable{X}$},
called the \emph{transfer of $f$}. (Actually, a more general
transfer for fibrations with compact fibres has been constructed
by Dwyer in \cite{Dw:Transfer} but we need not consider that
here.) We will use the following important properties of
transfers:

\medskip
{\bf Contravariant functoriality}\qua If \mbox{$f \negmedspace : X
\to Y$} and \mbox{$g \negmedspace : Y \to Z$} are finite coverings
of connected spaces, then
$$tr_{g \circ f} \simeq tr_f \circ tr_g.$$

{\bf Normalization}\qua If \mbox{$f \negmedspace : X \to Y$} is an
$n$--fold cover of a connected space, then the induced  map in
singular cohomology (with any coefficients)
  $$tr_f^* \circ \PtdStable{f}^* : \Coh{Y} \to \Coh{X} \to \Coh{Y}$$
is multiplication by $n$.

\medskip
{\bf Frobenius reciprocity}\qua If \mbox{$f \negmedspace : X
\to Y$} is a finite cover of a connected space, then the following
diagram, where $\Delta_X$ and $\Delta_Y$ denote the respective
diagonals of $X$ and $Y$, commutes:
\[
\begin{CD}
{\PtdStable{Y}} @ > {\PtdStable\Delta_Y} >> {\PtdStable{Y}\wedge\PtdStable{Y}} \\
@ VV {tr_f} V @ VV {1 \wedge tr_f} V \\
{\PtdStable{X}} @> {(\PtdStable{f}\wedge id) \circ \PtdStable\Delta_X} >> {\PtdStable{Y}\wedge\PtdStable{X}.}\\
\end{CD}
\]
In particular,
$$tr_f^*( f^*(y)\cdot x ) = y \cdot tr_f^*(x) $$
for \mbox{$x \in \Coh{X}$} and \mbox{$y \in \Coh{Y}.$}

\medskip
Since \mbox{$\PtdStable{X} \simeq \Stable{X} \vee \SphereSpectrum$},
the transfer \mbox{$tr_f \negmedspace : \PtdStable{Y} \to
\PtdStable{X}$} of a finite cover\break \mbox{$f \negmedspace : X \to  Y$}
restricts to a \emph{reduced transfer} \mbox{$\Stable{X} \to
\Stable{Y}$}. As there is no danger of confusion we also denote the
reduced transfer by $tr_f$, and sometimes refer to it as transfer.

A monomorphism of groups \mbox{$\psi \negmedspace : G \to G'$}
induces a fibration $G'/\psi(G) \hookrightarrow BG
\xrightarrow{B\psi} BG' $. If $[G':\psi(G)]$ is finite, which is
always the case if $G$ and $G'$ are finite, $B\psi$ therefore
admits a transfer map, which we denote $tr_{\psi}$ for short. In
the special case of the inclusion \mbox{$H \leq G$} of a subgroup
of finite index we denote the transfer by $tr_H$.

Given a finite $(G_1,G_2)$--biset \mbox{$\Omega \in
\Morita(G_1,G_2)$}, we now get a stable map\break
  \mbox{$\alpha(\Omega) \in \PtdStableMaps{BG_1}{BG_2}$}
as follows. Let \mbox{$\Lambda := G_2 \backslash \Omega$}. Since
the left $G_2$--action on $\Omega$ is free, we get a principal
fibre sequence
$$G_2 \to\Omega \times_{G_1} EG_1 \to \Lambda \times_{G_1} EG_1 .$$
$$\xi \negmedspace:  \Lambda \times_{G_1} EG_1 \to BG_2\leqno{\hbox{Let}}$$
be the classifying map of this fibration. The projection map
$$\Lambda \times_{G_1} EG_1 \to BG_1$$
is a finite covering. Let \mbox{$\tau \negmedspace: \PtdStable{BG_1}
\to \PtdStable{(\Lambda \times_{G_1} EG_1})$} be the associated
transfer map. The map $\alpha(\Omega)$ is now defined as
$$\alpha(\Omega) := \PtdStable{\xi} \circ \tau.$$
This assignment extends to a homomorphism
$$\alpha \negmedspace: A(G_1,G_2) \to \PtdStableMaps{{BG_1}}{{BG_2}}$$
of abelian groups. Although it may not be immediate from the
definition, it is shown for example in \cite{BF} that the map
$\alpha$ is natural in the sense that it sends the pairing of
\ref{eq:BurnsidePairing} to the composition pairing for stable
maps:
$$\alpha(\Omega \circ \Lambda) = \alpha(\Omega) \circ \alpha(\Lambda).$$
Thus $\alpha$ is a ring homomorphism when \mbox{$G_1=G_2$}. One can
check that the value of $\alpha$ on a basis element
\mbox{$[H,\psi]$} is
$$\alpha([H,\psi]) = \PtdStable{B\psi} \circ tr_{H}.$$
The homomorphism $\alpha$ gives a way to relate $A(G_1,G_2)$ to
the group of homotopy classes of stable maps
$\PtdStableMaps{{BG_1}}{{BG_2}}$. Lewis--May--McClure have made this
relationship precise in \cite{LMM}. As a consequence of the Segal
conjecture (proved by Carlsson in \cite{Car}), they show that
$\alpha$ is completion with respect to the augmentation ideal
$I(G_1)$ of the Burnside ring $A(G_1)$. In the case where $G_1$ is
a $p$--group, May--McClure \cite{MM} showed that, after getting rid
of basepoints, this completion takes a simple form, which we will
describe below.

\begin{definition}
For finite groups $G_1$ and $G_2$, we say that a $(G_1,G_2)$--pair
$(H,\psi)$ is \emph{trivial} if $\psi$ is the trivial
homomorphism. In this case we also say that the conjugacy class
$[H,\psi]$ is trivial. When $\psi$ is not the trivial
homomorphism, we say that the pair $(H,\psi)$ and the conjugacy
class $[H,\psi]$ are \emph{non-trivial}.
\end{definition}

Let $\wtilde{A}(G_1,G_2)$ be the quotient module obtained from
$A(G_1,G_2)$ by quotienting out all trivial basis elements
$[H,\psi]$. Recalling that \mbox{$\PtdStable{BG} \simeq
\Stable{BG}\vee \SphereSpectrum$}, where\break \mbox{$\SphereSpectrum =
\Stable{S^0}$} is the suspension sphere spectrum, one can check that
there is an induced map
$$\alpha \negmedspace :  \wtilde{A}(G_1,G_2) \longrightarrow \PtdStableMaps{BG_1}{BG_2} / \StableMaps{{BG_1}_+}{S^0} \cong \StableMaps{BG_1}{BG_2}.$$
May--McClure proved that when $G_1$ is a $p$--group, $I(G_1)$--adic
completion coincides with $p$--adic completion on
$\wtilde{A}(G_1,G_2)$, and deduced the following version of the
Segal conjecture.
\begin{theorem}[Segal conjecture \cite{Car,LMM,MM}] \label{thm:Segal}
If $S$ is a finite $p$--group and $G$ any finite group, then
$\alpha$ induces an isomorphism
$$\wtilde{\alpha} \negmedspace : \pComp{\wtilde{A}(S,G)} \stackrel{\cong}{\longrightarrow} \StableMaps{BS}{BG},$$
where \mbox{$\pComp{(-)} =(-)\otimes\Zp$} is $p$--adic completion.
\end{theorem}

For finite $p$--groups $S_1$ and $S_2$, we will, in view of the
Segal conjecture, have reason to $p$--complete the Burnside module
$A(S_1,S_2)$. The resulting $\Zp$--module $\pComp{A(S_1,S_2)}$ is a
free $\Zp$--module with one basis element for each conjugacy class
of $(S_1,S_2)$--pairs, and by a further, yet slight, abuse of
notation, we will also let $[P,\psi]_{S_1}^{S_2}$ (or $[P,\psi]$)
denote the basis element of $\pComp{A(S_1,S_2)}$ corresponding to
the conjugacy class of the $(S_1,S_2)$--pair $(P,\psi)$.

We conclude this section by adapting some ``bookkeeping'' tools for
$(S_1,S_2)$--bisets to keep track of elements of
$\pComp{A(S_1,S_2)}$. First we note that the structure of
$\pComp{A(S_1,S_2)}$ allows us to define a collection of
homomorphisms
$$\chi_{[P,\psi]} \negmedspace : \pComp{A(S_1,S_2)} \to \Zp,$$
one for each conjugacy class of $(S_1,S_2)$--pairs, by demanding
that
$$\Omega = \sum_{[P,\psi]} \chi_{[P,\psi]}(\Omega)\cdot[P,\psi],$$
for all \mbox{$\Omega \in \pComp{A(S_1,S_2)}$}.

Next, we extend the notion of counting the number of $S_2$--orbits of
$(S_1,S_2)$--bisets to obtain a form of augmentation for Burnside
modules. The resulting assignment is natural in that it sends the
pairing of \ref{eq:BurnsidePairing} to multiplication in $\Zp$.

\begin{lemma}\label{lem:StableAugmentation}
For every pair of finite $p$--groups $S_1$ and $S_2$, the
assignment
$$\Morita(S_1,S_2) \to \mathbb{Z},\hspace{.3cm} \Omega \mapsto |S_2 \backslash \Omega|$$
extends to a homomorphism
$$\countS\negmedspace : \pComp{A(S_1,S_2)} \to \Zp,$$
sending composition to multiplication.
\end{lemma}
\begin{proof}
Recalling that bisets \mbox{$\Omega \in \Morita(S_1,S_2)$} have a
free $S_2$--action, we see that each assignment
$$\Morita(S_1,S_2) \to \mathbb{Z},\hspace{.3cm} \Omega \mapsto |S_2 \backslash \Omega|$$
is a morphism of monoids, and so we get an induced homomorphism
\mbox{$\countS\negmedspace : A(S_1,S_2) \to \mathbb{Z}$} and,
after $p$--completion, an induced homomorphism
\mbox{$\countS\negmedspace : \pComp{A(S_1,S_2)} \to \Zp$}.

Using the freeness of the left action for bisets \mbox{$\Omega \in
\Morita(S_2,S_3)$} and \mbox{$\Lambda \in \Morita(S_1,S_2)$}
again, we get
$$|S_3 \backslash (\Omega \circ \Lambda)| = |\Omega \times_{S_2} \Lambda|/|S_3| = (|\Omega|\cdot|\Lambda|/|S_2|)/|S_3| = |S_3 \backslash \Omega| \cdot |S_2 \backslash \Lambda|.$$
The collection of homomorphisms \mbox{$\countS\negmedspace :
\pComp{A(S_1,S_2)} \to \Zp$} therefore sends composition to
multiplication.
\end{proof}

A useful, well known result states that for a finite group $G$, two
finite $G$--sets $\Omega$ and $\Lambda$ are isomorphic if and only if
they have the same number of fixed points for every subgroup of $G$.
Since the number of fixed points depends only on the conjugacy class
of the subgroup, an alternative formulation is that there is an
injective $\mathbb{Z}$--module homomorphism
$$A(G) \to \prod_{[H]} \mathbb{Z}, \hspace{0.3cm} \Omega \mapsto \prod_{[H]} |\Omega^H|,$$
where the product is taken over conjugacy classes of subgroups
\mbox{$H \leq G$}.

For finite groups $S_1$ and $S_2$ we regard a $(S_1,S_2)$--biset
$\Omega$ as a left \mbox{$(S_1 \times S_2)$}--set by putting
\mbox{$(g,h)x := hxg^{-1}$} for \mbox{$g \in S_1$}, \mbox{$h \in
S_2$} and \mbox{$x \in \Omega$}. This assignment preserves
isomorphism classes and we obtain an injection
$$A(S_1,S_2) \longrightarrow A(S_1 \times S_2)$$
sending a basis element $[P,\psi]$ to \mbox{$[(S_1 \times
S_2)/\Graph{P}{\psi}].$} For a subgroup \mbox{$Q \leq S_1 \times
S_2$}, this allows us to define $\Omega^Q$ as the fixed-point set of
$\Omega$ under the action of $Q$. By linear extension and
$p$--completion we get a well defined $\Zp$--module homomorphism
$$\pComp{A(S_1,S_2)} \longrightarrow \Zp, \hspace{.3cm} \Omega \mapsto |\Omega^Q|,$$
depending only on the conjugacy class of $Q$. On basis elements we
have
$$\left|[P,\psi]^{Q}\right| = \left|\left((S_1\times S_2)/\Graph{P}{\psi}\right)^Q \right| = \left|\Graph{P}{\psi}\backslash N_{S_1\times S_2}(Q,\Graph{P}{\psi})\right| = \frac{\left| N_{S_1\times S_2}(Q,\Graph{P}{\psi}) \right|}{\left| P \right| }.$$

\begin{lemma}\label{lem:StableFxdPts}
Let $S_1$ and $S_2$ be finite $p$--groups. Then the $\Zp$--module
homomorphism
$$\pComp{A(S_1,S_2)} \longrightarrow \prod_{[P,\psi]} \Zp, \hspace{.3cm} \Omega \mapsto \prod_{[P,\psi]} \big|\Omega^{\Graph{P}{\psi}}\big|,$$
where the product is taken over conjugacy classes of
$(S_1,S_2)$--pairs, is injective.
\end{lemma}
\begin{proof}
Being a composition of two injective homomorphisms, the
$\Z$--module homomorphism
$$A(S_1,S_2) \longrightarrow A(S_1\times S_2) \longrightarrow \prod_{[Q]} \mathbb{Z}, \hspace{1cm} \Omega \mapsto \prod_{[Q]} \left|\Omega^Q\right|,$$
where the product runs over conjugacy classes of subgroups
\mbox{$Q \leq S_1 \times S_2$}, is itself injective. Noting that
the collection of graphs of $(S_1,S_2)$--pairs is closed under
conjugation in \mbox{$S_1\times S_2$} and taking subgroups, we see
that for a $(S_1,S_2)$--pair $(P,\psi)$, we have
$$N_{S_1\times S_2}(Q,\Graph{P}{\psi}) = \emptyset,$$
and consequently
$$\left|[P,\psi]^{Q}\right| = 0,$$
if $Q$ is not the graph of an $(S_1,S_2)$--pair. We conclude that the
restriction to a $\Z$--module homomorphism
$$A(S_1,S_2) \longrightarrow \prod_{[P,\psi]} \mathbb{Z}, \hspace{1cm} \Omega \mapsto \prod_{[P,\psi]} \big|\Omega^{\Graph{P}{\psi}}\big|,$$
is injective, and it remains so after $p$--completion.
\end{proof}

\section{Fusion subconjugacy} \label{sec:Subconjugacy}
In this section we introduce the notion of fusion subconjugacy for
subgroups of a finite $p$--group $S$ and for $(S_1,S_2)$--pairs. This
induces a filtration on the $p$--completed Burnside module
$\pComp{A(S_1,S_2)}$ and consequently of the group
$\StableMaps{BS_1}{BS_2}$ of homotopy classes of stable maps. By
studying this filtration we will obtain useful information about how
homotopy classes of stable maps between classifying spaces of finite
\mbox{$p$--groups} behave under composition with stable maps arising
from fusion systems over those \mbox{$p$--groups}. The material in
this section is presented for $p$--completed Burnside modules because
we are mostly interested in that setting. However the analogous
results still hold in the uncompleted or $p$--localized case.

\begin{definition} \label{def:GroupSubconjugacy}
Let $\F$ be a fusion system over  a finite $p$--group $S$, and let
$P$ and $Q$ be subgroups of $S$.
\begin{itemize}
  \item
We say that $Q$ is $\F$--\emph{subconjugate} to $P$, and write
\mbox{$Q \Fsub{\F} P$}, if there exists a morphism \mbox{$\varphi
\in \HomF{Q}{P}$}.
  \item
We say that $Q$ is $\F$--\emph{conjugate} to $P$, and write
\mbox{$Q \Fcon{\F} P$} if $Q$ is isomorphic to $P$ in $\F$.
  \item
We say that $Q$ is \emph{strictly} \mbox{$\F$--\emph{subconjugate}}
to $P$, and write \mbox{$Q \Fstrictsub{\F} P$} if $Q$ is
\mbox{$\F$--subconjugate} to $P$, but not \mbox{$\F$--conjugate} to
$P$.
\end{itemize}
\end{definition}
When there is no danger of confusion, we will write $\noFsub,
\noFstrictsub,$ and $\noFcon$ instead of $\Fsub{\F},
\Fstrictsub{\F}$ and $\Fcon{\F}$. For the fusion system $\F_S$ of
$S$, $\F_S$--subconjugacy coincides with $S$--subconjugacy.

We make a similar definition for pairs.
\begin{definition} \label{def:PairSubconjugacy}
Let $\F_1$ and $\F_2$ be fusion systems over finite $p$--groups
$S_1$ and $S_2$, respectively. Let $(P,\psi)$ and $(Q,\rho)$ be
two $(S_1,S_2)$--pairs.
\begin{itemize}
  \item
We say that $(Q,\rho)$ is $(\F_1,\F_2)$--\emph{subconjugate} to
$(P,\psi)$, and write
$$(Q,\rho) \Fsub{(\F_1,\F_2)} (P,\psi),$$
if there exist morphisms \mbox{$\varphi_1 \in \Hom{Q}{P}{\F_1}$}
and \mbox{$\varphi_2 \in \Hom{\rho(Q)}{\psi(P)}{\F_2}$} such that
the following diagram commutes
\[
\begin{CD}
{Q} @ > {\rho} >> {\rho(Q)} \\
@ V \varphi_1 VV @ VV \varphi_2 V \\
{P} @> {\psi} >> {\psi(P).} \\
\end{CD}
\]
  \item
We say that $(Q,\rho)$ is $(\F_1,\F_2)$--\emph{conjugate} to
$(P,\psi)$, and write
$$(Q,\rho) \Fcon{(\F_1,\F_2)} (P,\psi),$$
if
$$(Q,\rho) \Fsub{(\F_1,\F_2)} (P,\psi) \hspace{1cm} \text{and} \hspace{1cm} (P,\psi) \Fsub{(\F_1,\F_2)} (Q,\rho).$$
  \item
We say that $(Q,\rho)$ is \emph{strictly}
$(\F_1,\F_2)$--\emph{subconjugate} to $(P,\psi)$, and write
$$(Q,\rho)
\Fstrictsub{(\F_1,\F_2)} (P,\psi),$$ if $(Q,\rho)$ is
$(\F_1,\F_2)$--subconjugate to $(P,\psi)$, but not
\mbox{$(\F_1,\F_2)$--conjugate} to $(P,\psi)$.
\end{itemize}
\end{definition}
When there is no danger of confusion, we will write $\noFsub$,
$\noFstrictsub,$ and $\noFcon$ instead of\break $\Fsub{(\F_1,\F_2)}$,
$\Fstrictsub{(\F_1,\F_2)}$ and $\Fcon{(\F_1,\F_2)}$. As before,
$(\F_{S_1},\F_{S_2})$--conjugacy agrees with the notion of
$(S_1,S_2)$--conjugacy defined in \fullref{sec:Burnside}.

\begin{remark} \label{rem:SRespectSubconjugacy}
It is easy to check that $(\F_1,\F_2)$--subconjugacy is preserved by
$(S_1,S_2)$--conjugacy. Therefore we will often say that an
isomorphism class of pairs $[Q,\rho]$ is $(\F_1,\F_2)$--subconjugate
to an isomorphism class $[P,\psi]$ and write
$$[Q,\rho] \Fsub{(\F_1,\F_2)} [P,\psi] \hspace{0.5cm} \text{(or}~ [Q,\rho] \noFsub [P,\psi] \text{)}  $$
if the subconjugacy relation
$$(Q,\rho) \Fsub{(\F_1,\F_2)} (P,\psi)$$
holds between any (and hence all) representatives of the classes.
Furthermore, we will use the same terminology when we regard
$[P,\psi]$ and $[Q,\rho]$ as basis elements of $A(S_1,S_2)$ or
$\pComp{A(S_1,S_2)}$. The analogous remark applies to
$(\F_1,\F_2)$--conjugacy and strict $(\F_1,\F_2)$--subconjugacy.
\end{remark}

\begin{remark} \label{rem:SubconjugacyLink}
Subconjugacy among $(S_1,S_2)$--pairs can in fact be regarded as a
special case of subconjugacy among subgroups of \mbox{$S_1 \times
S_2$}. Recall from \cite[Section 1]{BLO2}, that in the setting of
the definition above, the fusion system \mbox{$\F_1 \times \F_2$}
over \mbox{$S_1 \times S_2$} is defined by setting
\mbox{$\Hom{P}{Q}{\F_1 \times \F_2}$} to be the morphism set
%$$\Hom{P}{Q}{\F_1 \times \F_2} := \{(\varphi_1,\varphi_2)\vert_P \in \Hom{P}{Q}{} \mid \varphi_i \in \Hom{P_i}{S_i}{\F_i}, P\leq P_1 \times P_2 \}$$
$$\{(\varphi_1,\varphi_2)\vert_P \in \Hom{P}{Q}{} \mid \varphi_i \in \Hom{P_i}{S_i}{\F_i}, P\leq P_1 \times P_2 \}$$
for all \mbox{$P,Q \leq S_1 \times S_2$}. By \cite[Lemma 1.5]{BLO2},
the fusion system \mbox{$\F_1 \times \F_2$} is saturated if the
fusion systems $\F_1$ and $\F_2$ are both saturated. For
$(S_1,S_2)$--pairs $(P,\psi)$ and $(P',\psi')$, one can check that
$(P',\psi')$ is \mbox{$(\F_1,\F_2)$--subconjugate} to $(P,\psi)$ if
and only $\Graph{P'}{\psi'}$ is $(\F_1 \times \F_2)$--subconjugate to
$\Graph{P}{\psi}$.
\end{remark}

It is easy to check that $(\F_1,\F_2)$--subconjugacy is a transitive
relation. Therefore the $(\F_1,\F_2)$--conjugacy classes of
$(S_1,S_2)$--pairs form a poset under $(\F_1,\F_2)$--subcon\-jugacy.
Since $(S_1,S_2)$--conjugacy classes of $(S_1,S_2)$--pairs form a
$\Zp$--basis for the $\Zp$--module $\pComp{A(S_1,S_2)}$, this leads us
to a poset-indexed filtration as defined below.
\begin{definition} \label{def:SubconjugacyFiltration}
Let $\F_1$ and $\F_2$ be fusion systems over finite $p$--groups
$S_1$ and $S_2$, respectively. Let $(P,\psi)$ be a
$(S_1,S_2)$--pair.
\begin{itemize}
  \item
Let $\FsubModule{\F_1}{\F_2}{P}{\psi}$ denote the submodule of
$\pComp{A(S_1,S_2)}$ generated by the basis elements $[Q,\rho]$
such that
$$[Q,\rho] \Fsub{(\F_1,\F_2)} [P,\psi].$$
  \item
Let $\FconModule{\F_1}{\F_2}{P}{\psi}$ denote the submodule of
$\pComp{A(S_1,S_2)}$ generated by the basis elements $[Q,\rho]$
such that
$$[Q,\rho] \Fcon{(\F_1,\F_2)} [P,\psi].$$
  \item
Let $\FstrictsubModule{\F_1}{\F_2}{P}{\psi}$ denote the submodule
of $\pComp{A(S_1,S_2)}$ generated by the basis elements $[Q,\rho]$
such that
$$[Q,\rho] \Fstrictsub{(\F_1,\F_2)} [P,\psi].$$
\end{itemize}
\end{definition}
When the fusion systems $\F_1$ and $\F_2$ are clear from the
context, and there is no danger of confusion, we will write
$\noFsubModule{P}{\psi}, \noFconModule{P}{\psi}$ and
$\noFstrictsubModule{P}{\psi}$ instead of
$\FsubModule{\F_1}{\F_2}{P}{\psi}, \FconModule{\F_1}{\F_2}{P}{\psi}$
and $\FstrictsubModule{\F_1}{\F_2}{P}{\psi}$.

The stable selfmaps of a finite $p$--group arising from morphisms in
a fusion system are of special importance in this paper. We
therefore consider the corresponding subring of $\pComp{A(S,S)}$.
\begin{definition}
Let $\F$ be a fusion system over a finite $p$--group $S$. We denote
by $A_{\F}(S,S)$ the submodule of $A(S,S)$ generated by the basis
elements $[P,\varphi]$ where \mbox{$\varphi \in \HomF{P}{S}$}.
\end{definition}
After $p$--completion we obtain a submodule $\StableFusion{\F}{S}$
of $\pComp{A(S,S)}$, again generated by the basis elements
$[P,\varphi]$ where \mbox{$\varphi \in \HomF{P}{S}$}.

\begin{remark} \label{rem:RFDescription}
One can check that
$$  \StableFusion{\F}{S} = \FsubModule{\F_S}{\F}{S}{id} = \FsubModule{\F}{\F}{S}{id} = \FsubModule{\F}{\F_S}{S}{id}.$$
\end{remark}

Under composition, the $\Zp$--module $\pComp{A(S_1,S_2)}$ becomes a
left $\pComp{A(S_2,S_2)}$--module and a right
$\pComp{A(S_1,S_1)}$--module. The filtration in 
\fullref{def:SubconjugacyFiltration} is useful to us mainly because of
the following lemma.
\begin{lemma} \label{lem:PreservesFiltration}
Let $\F_1$ and $\F_2$ be fusion systems over the finite $p$--groups
$S_1$ and $S_2$, respectively. The following hold for every
$(S_1,S_2)$--pair $(P,\psi)$:
\begin{enumerate}
  \item[(a)] $\StableFusion{\F_2}{S_2} \circ\noFsubModule{P}{\psi} \subseteq  \noFsubModule{P}{\psi},$
  \item[(b)] $\StableFusion{\F_2}{S_2} \circ \noFstrictsubModule{P}{\psi} \subseteq \noFstrictsubModule{P}{\psi},$
  \item[(c)] $\noFsubModule{P}{\psi} \circ \StableFusion{\F_1}{S_1} \subseteq \noFsubModule{P}{\psi},$
  \item[(d)] $\noFstrictsubModule{P}{\psi} \circ \StableFusion{\F_1}{S_1} \subseteq \noFstrictsubModule{P}{\psi}.$
\end{enumerate}
\end{lemma}
\begin{proof}
We prove parts (a) and (b), and leave the proofs of (c) and (d),
which are similar, to the reader.

First we show that for any $(S_1,S_2)$--pair $(Q,\rho)$ and any basis
element $[T,\varphi]$ of $\StableFusion{\F_2}{S_2}$, we have
$$[T,\varphi] \circ [Q,\rho] \in \noFsubModule{Q}{\rho}.$$
Indeed, by the double coset formula,
$$[T,\varphi] \circ [Q,\rho] = \DCF{T}{\varphi}{Q}{\rho}{S_2}, $$
and it suffices to prove that
$$\DCFterm{T}{\varphi}{Q}{\rho} \noFsub [Q,\rho]$$
for each \mbox{$x\in S_2$}. But this is clear by the diagram
\[
\begin{CD}
{\rho^{-1}\left(\rho\left(Q\right)\cap T^x\right)} @ > \varphi\circ c_x\circ{\rho} >> {\varphi\left(\lsup{\rho\left(Q\right)}{x} \cap T\right)} \\
@ V \iota VV @ VV c_x^{-1} \circ \varphi^{-1} V  \\
{Q} @ > {\rho} >> {\rho\left(Q\right).} \\
\end{CD}
\]
To prove part (a), let \mbox{$[Q,\rho] \noFsub [P,\psi]$}. By the
preceding observation we get
$$[T,\varphi] \circ [Q,\rho] \in \noFsubModule{Q}{\rho} \subseteq \noFsubModule{P}{\psi}.$$
Letting $[T,\varphi]$ and $[Q,\rho]$ vary over all basis elements of
$\StableFusion{\F_2}{S_2}$ and\break $\noFsubModule{P}{\psi}$, we get the
desired result.

Similarly, part (b) follows upon noting that for a basis element
$[Q,\rho]$ of $\noFstrictsubModule{P}{\psi}$ and a basis element
$[T,\varphi]$ of $\StableFusion{\F_2}{S_2}$, we have
$$[T,\varphi] \circ [Q,\rho] \in \noFsubModule{Q}{\rho} \subseteq \noFstrictsubModule{P}{\psi}.\proved$$
\end{proof}

We have the following structural corollaries.
\begin{corollary} \label{cor:RFclosed}
Let $\F$ be a fusion system over a finite $p$--group $S$. Then
$\StableFusion{\F}{S}$ is a subring of $\pComp{A(S,S)}$. Similarly
$A_{\F}(S,S)$ is a subring of $A(S,S)$.
\end{corollary}
\begin{proof}
This follows from \fullref{rem:RFDescription} and 
\fullref{lem:PreservesFiltration}. The same proof works for the last
statement.
\end{proof}

\begin{corollary}
Let $\F_1$ and $\F_2$ be fusion systems over the finite $p$--groups
$S_1$ and $S_2$, respectively. For every $(S_1,S_2)$--pair
$(P,\psi)$, the $\Zp$--modules $\noFsubModule{P}{\psi}$ and
$\noFstrictsubModule{P}{\psi}$ are left modules over
$\StableFusion{\F_2}{S_2}$ and right modules over
$\StableFusion{\F_1}{S_1}$.
\end{corollary}

\begin{definition}
Let $\F$ be a fusion system over a finite $p$--group $S$, and let
\mbox{$\Omega \in \pComp{A(S,S)}$}. We say that $\Omega$ is
\emph{right $\F$--stable} if for every \mbox{$P \leq S$} and every
\mbox{$\varphi \in \HomF{P}{S}$} we have
$$\Omega \circ [P,\varphi]_P^{S} = \Omega \circ [P,\iota_P]_P^{S}$$
in $\pComp{A(P,S)}$. Similarly we say that $\Omega$ is \emph{left
$\F$--stable} if for every \mbox{$P \leq S$} and \mbox{$\varphi \in
\HomF{P}{S}$} we have
$$[\varphi(P),\varphi^{-1}]_{S}^P \circ \Omega = [P,id_P]_{S}^{P} \circ \Omega$$
in $\pComp{A(S,P)}$.
\end{definition}
When $\Omega$ is represented by a $(S,S)$--biset $X$, the right
$\F$--stability condition means that the restriction of $X$ to a
$(P,S)$--biset via $\varphi$ is isomorphic to the restriction of $X$
via the inclusion \mbox{$P \hookrightarrow S$}, while left stability
means that the restriction of $X$ to a $(S,P)$--biset via $\varphi$
is isomorphic to the restriction of $X$ via the inclusion.

We will later define a similar notion of fusion stability for maps
between stable classifying spaces of $p$--groups.

\begin{lemma} \label{lem:FinvOnConjugates}
Let $\F_1$ and $\F_2$ be fusion systems over the finite $p$--groups
$S_1$ and $S_2$, respectively, let \mbox{$\Omega_1 \in
\pComp{A(S_1,S_1)}$} be left $\F_1$--stable, and let \mbox{$\Omega_2
\in \pComp{A(S_2,S_2)}$} be right $\F_2$--stable. If the
$(S_1,S_2)$--pairs $(P,\psi)$ and $(Q,\rho)$ are
$(\F_1,\F_2)$--conjugate, then
$$\Omega_2 \circ [Q,\rho] \circ \Omega_1 = \Omega_2 \circ [P,\psi] \circ \Omega_1.$$
\end{lemma}
\begin{proof}
Let \mbox{$\wtilde{\psi} \negmedspace : P \to \psi(P)$} denote
the restriction of $\psi$ to its image. Since $(P,\psi)$ and
$(Q,\rho)$ are conjugate, there exist isomorphisms
\mbox{$\varphi_1 \in \Hom{P}{Q}{\F_1}$} and \mbox{$\varphi_2 \in
\Hom{\psi(P)}{\rho(Q)}{\F_2}$} such that
$$\rho \circ \varphi_1  = \iota_{\rho(Q)} \circ \varphi_2 \circ \wtilde{\psi}.$$
Using stability, and recalling the simple description of the double
coset formula for composable morphisms in \ref{eq:DCFsimple}, we
now obtain
\begin{align*}
  &\Omega_2 \circ [Q,\rho]_{S_1}^{S_2} \circ \Omega_1 \\
  &\hspace{13pt} = \Omega_2 \circ [Q,\iota_{\rho(Q)} \circ \varphi_2\circ \wtilde{\psi} \circ \varphi_1^{-1}]_{S_1}^{S_2} \circ \Omega_1 \\
  &\hspace{13pt} = \Omega_2 \circ \left( [\psi\left(P\right),\iota_{\rho(Q)} \circ \varphi_2]_{\psi\left(P\right)}^{S_2} \circ [P,\wtilde{\psi}]_{P}^{\psi(P)} \circ [Q,\varphi_1^{-1}]_{S_1}^{P}\right) \circ \Omega_1 \\
  &\hspace{13pt} = \left(\Omega_2 \circ [\psi\left(P\right),\iota_{\rho(Q)} \circ \varphi_2]_{\psi\left(P\right)}^{S_2}\right) \circ [P,\wtilde{\psi}]_{P}^{\psi\left(P\right)} \circ \left([\varphi(P),\varphi_1^{-1}]_{S_1}^{P} \circ \Omega_1 \vphantom{{]}^{()}}\right) \\
  &\hspace{13pt} = \left(\Omega_2 \circ [\psi\left(P\right),\iota_{\psi\left(P\right)}]_{\psi\left(P\right)}^{S_2}\right) \circ [P,\wtilde{\psi}]_{P}^{\psi\left(P\right)} \circ \left([P,id_P]_{S_1}^{P} \circ \Omega_1 \vphantom{{]}^{()}}\right) \\
  &\hspace{13pt} = \Omega_2 \circ \left([\psi\left(P\right),\iota_{\psi\left(P\right)}]_{\psi\left(P\right)}^{S_2} \circ [P,\wtilde{\psi}]_{P}^{\psi\left(P\right)} \circ [P,id_P]_{S_1}^{P}) \right) \circ \Omega_1 \\
  &\hspace{13pt} = \Omega_2 \circ [P,\psi]_{S_1}^{S_2} \circ \Omega_1.
\end{align*}
This completes the proof. \end{proof}

\section{Characteristic idempotents} \label{sec:Idempotent}
In this section, and for the rest of the paper, we restrict our
attention to saturated fusion systems. For a saturated fusion system
$\F$ over a finite \mbox{$p$--group $S$}, we will prove the existence
of an idempotent \mbox{$\CharIdem \in \pComp{A(S,S)}$}, related to
$\F$ through properties made precise in 
\fullref{def:CharIdempotent} below. These properties, and their
importance, were originally recognized by Linckelmann--Webb for
bisets. It is the careful analysis of $\CharIdem$ which will allow
us to produce the main results of this paper. In later sections we
will see that $\CharIdem$ is uniquely determined by $\F$ and that it
characterizes the fusion system $\F$, thus justifying the term
characteristic idempotent.

In \cite[Section 5]{BLO2}, Broto--Levi--Oliver determined the
cohomological structure of a $p$--local finite group $(S,\F,\Link)$.
In short, they proved that in cohomology, the natural inclusion
\mbox{$\theta \negmedspace : BS \to \ClSp$} induces an isomorphism
$$\Coh{\ClSp} \stackrel{\cong}{\longrightarrow} \Coh{\F} \subseteq{\Coh{BS}},$$
where
$$
 \Coh{\F} := \invlimlim{\F}\Coh{B(-)}
$$
is the ``ring of stable elements for \F'', regarded as a subring of
$\Coh{BS}$, via the identification
$$
  \Coh{\F} \cong \{ x \in \Coh{BS} \mid B\varphi^*(x) = B\iota_P^*(x)~ \text{for all}~ P \leq S, \varphi \in \HomF{P}{S}   \}.
$$
One of the key ingredients in their proof is the construction of a
characteristic biset \mbox{$\Omega \in \Morita(S,S)$}, as defined
below. We take advantage of their construction and produce our
characteristic idempotent by showing the convergence of a
judiciously chosen subsequence of the sequence
$$[\Omega], [\Omega]^2, [\Omega]^3, \dots$$

\begin{definition} Let $\F$ be a fusion system over a finite
$p$--group $S$. We say that an element \mbox{$\Omega \in A(S,S)$} is
a virtual characteristic biset for $\F$ if it has the following
properties:
\begin{enumerate}
  \item[$\rm(a')$] $\Omega \in A_{\F}(S,S)$.
  \item[$\rm(b1')$] $\Omega$ is right $\F$--stable.
  \item[$\rm(b2')$] $\Omega$ is left $\F$--stable.
  \item[$\rm(c')$] $\countS(\Omega) \equiv 1~ (mod~ p)$.
\end{enumerate}
If in addition \mbox{$\Omega \in \Morita(S,S)$} then we say that
$\Omega$ is a characteristic biset for $\F$.
\end{definition}
We refer to these properties as the Linckelmann--Webb properties as
they were first suggested in unpublished work of Linckelmann--Webb
\cite{LW}, although Property $\rm(b2')$ did not feature there. We refer
to Properties $\rm(b1')$ and $\rm(b2')$ collectively as Property $\rm(b')$. The
Linckelmann--Webb properties mimic the properties of a finite group
$G$ with Sylow subgroup $S$ regarded as an $(S,S)$--biset, although
some scaling may be required to obtain Property $\rm(c')$. The importance
of the Linckelmann--Webb properties is apparent in the following
result.

\begin{proposition}{\rm\cite{LW,BLO2}}\qua \label{prop:LWImportance}
Let $\F$ be a fusion system over a finite $p$--group $S$. If $\Omega$
is a virtual characteristic biset for $\F$, then the induced map
$\alpha(\Omega)^*$ in cohomology is an idempotent in
$End(\Coh{BS})$, is $\Coh{\F}$--linear and a homomorphism of modules
over the Steenrod algebra; and
$$Im[\Coh{BS}\xrightarrow{\alpha(\Omega)^*}\Coh{BS}] =
\Coh{\F}.$$
\end{proposition}
\begin{proof} See the proof of \cite[Proposition 5.5]{BLO2}.
\end{proof}

A characteristic idempotent for a fusion system $\F$ over $S$ is an
idempotent in $\pComp{A(S,S)}$ with $p$--completed, idempotent
analogues of the Linckelmann--Webb properties. This is stated
precisely below.
\begin{definition} \label{def:CharIdempotent}
Let $\F$ be a fusion system over a finite $p$--group $S$. A
\emph{characteristic idempotent for $\F$} is an idempotent
\mbox{$\CharIdem \in \pComp{A(S,S)}$} with the following
properties:
\begin{enumerate}
  \item[(a)] $\CharIdem \in \StableFusion{\F}{S}$.
  \item[(b1)] $\CharIdem$ is right $\F$--stable.
  \item[(b2)] $\CharIdem$ is left $\F$--stable.
  \item[(c)] $\countS(\CharIdem) = 1$.
\end{enumerate}
\end{definition}
We again refer to Properties (b1) and (b2) collectively as
Property (b).

The existence of characteristic bisets for saturated fusion systems
was established by Broto--Levi--Oliver in \cite{BLO2} through a
constructive argument. Although they, like Linckelmann--Webb, did not
include Property $\rm(b2')$ in their statement of the result, it is
implicit in their construction.
\begin{proposition}{\rm\cite[Proposition 5.5]{BLO2}}\qua\label{prop:BLOIdempotent}
Every saturated fusion system $\F$ over a $p$--group $S$ has a
characteristic $(S,S)$--biset.
\end{proposition}
The preceding proposition is the only point in this paper where we
rely on the saturation of fusion systems. If we were instead to
assume that every fusion system in sight has a characteristic
biset, then the construction of characteristic idempotents and
classifying spectra, as well as the proof of their properties
still go through. It is an interesting question whether this
really amounts to weakening our hypothesis. That is, whether the
existence of a characteristic biset for a fusion system $\F$
implies that $\F$ is saturated.

We now proceed by a sequence of lemmas about $(S,S)$--bisets to
produce the characteristic idempotent.
\begin{lemma} \label{lem:CanTakePowers}
Let $\Omega$ and $\Lambda$ be two (virtual) characteristic bisets
for a fusion system $\F$ over a finite $p$--group $S$. Then
\mbox{$\Omega \circ \Lambda$} is also a (virtual) characteristic
biset for $\F$. In particular, any power of $\Omega$ is a
(virtual) characteristic biset for $\F$.
\end{lemma}
\begin{proof}
That \mbox{$\Omega \circ \Lambda$} has Property $\rm(a')$ follows from
\fullref{cor:RFclosed}. To see that \mbox{$\Omega \circ
\Lambda$} has Property $\rm(b')$, we note that for \mbox{$P \leq S$} and
\mbox{$\varphi \in \HomF{P}{S}$} we have
$$
  (\Omega \circ \Lambda) \circ [P,\varphi]_P^S = \Omega \circ (\Lambda \circ [P,\varphi]_P^S) = \Omega \circ (\Lambda \circ [P,\iota_P]_P^S) = (\Omega \circ \Lambda) \circ
  [P,\iota_P]_P^S,
$$
and similarly
$$[\varphi(P),\varphi^{-1}]_S^P \circ (\Omega \circ \Lambda) = [P,id_P]_S^P \circ (\Omega \circ \Lambda).$$
Property $\rm(c')$ is clearly preserved since $\countS$ is
multiplicative. The final statement now follows by induction.
\end{proof}

\begin{lemma}\label{lem:ExistsM}
Let \mbox{$\Omega \in A(S,S)$}. Then there exists an $M>0$ such
that $\Omega^M$ is idempotent \mbox{$mod~ p$}.
\end{lemma}
\begin{proof}
Let $\bar{\Omega}$ denote the image of $\Omega$ under the
projection
$$A(S,S) \to A(S,S)/pA(S,S).$$
It is equivalent to show that $\bar{\Omega}^M$ is idempotent for
some \mbox{$M>0$}. Now, $A(S,S)$ is a finitely generated
$\mathbb{Z}$--module and hence $A(S,S)/pA(S,S)$ is finite. Consider
the sequence
$$\bar{\Omega}, \bar{\Omega}^2, \bar{\Omega}^3, \dots$$
in $A(S,S)/pA(S,S)$. By the pigeonhole principle there must be
numbers \mbox{$N,t>0$} such that \mbox{$\bar{\Omega}^N =
\bar{\Omega}^{N+t}$}. It follows that
$$\bar{\Omega}^n = \bar{\Omega}^{n+t}$$
for all \mbox{$n \geq N$}. Now take \mbox{$m \geq 0$} such that
$mt > N$ and put \mbox{$M := mt$}. Then
$$\bar{\Omega}^{2M} = \bar{\Omega}^{M+mt} = \bar{\Omega}^{M+(m-1)t} = \dots
=\bar{\Omega}^{M+t} = \bar{\Omega}^{M}.\proved$$
\end{proof}

The following two lemmas were suggested to the author by Bob
Oliver. Although they hold for any $p$--torsion-free ring, we will
state them only for $A(S,S)$.
\begin{lemma} \label{lem:IncreasingIdempotents}
If \mbox{$\Omega \in A(S,S)$} is idempotent \mbox{$mod~ p^k$},
where $k > 0$, then $\Omega^{p}$ is idempotent \mbox{$mod~
p^{k+1}$}.
\end{lemma}
\begin{proof}
Put $q:=p^k$. By assumption we can write
\begin{equation} \label{eq:IdemModp}
\Omega^2 = \Omega + q\Lambda
\end{equation}
for some \mbox{$\Lambda \in A(S,S)$}. It follows that
$$\Omega^2 + q\Omega\Lambda = \Omega(\Omega + q\Lambda) = \Omega^3
= (\Omega + q\Lambda)\Omega = \Omega^2 + q\Lambda\Omega,$$ so
$$q\Omega\Lambda = q\Lambda\Omega.$$
Since $A(S,S)$ is torsion-free as a $\mathbb{Z}$--module, we deduce
that $\Omega$ and $\Lambda$ commute. This allows us to apply the
binomial formula to \ref{eq:IdemModp} and get
$$\Omega^{2p} = \Omega^p + {{p}\choose{1}}\Omega^{p-1}q\Lambda + {{p}\choose{2}}\Omega^{p-2}q^2\Lambda^2 + \dots + {{p}\choose{p-1}}\Omega q^{p-1}\Lambda^{p-1} +q^p\Lambda^p.$$
A brief inspection of the coefficients occurring on the right hand
side, taking into account that $p$ divides $q$ since $k > 0$, shows
that we can therefore write
$$\Omega^{2p} = \Omega^p +pq\Lambda'$$
for some \mbox{$\Lambda' \in A(S,S)$}. Since \mbox{$pq = p^{k+1}$}
we deduce that $\Omega^{p}$ is idempotent \mbox{$mod~ p^{k+1}.$}
\end{proof}

\begin{lemma}\label{lem:Convergent}
If \mbox{$\Omega \in A(S,S)$} is idempotent \mbox{$mod~ p,$} then
the sequence
$$\Omega, \Omega^p, \Omega^{p^2}, \dots$$
converges in $\pComp{A(S,S)}$. Furthermore the limit is
idempotent.
\end{lemma}
\begin{proof}
By \fullref{lem:IncreasingIdempotents} and induction,
$\Omega^{p^k}$ is idempotent \mbox{$mod~ p^{k+1}$} for each
\mbox{$k \geq 0$}. That is to say,
\begin{equation} \label{eq:Convergent1}
\Omega^{2p^k} - \Omega^{p^k} \in p^{k+1}A(S,S)
\end{equation}
for \mbox{$k \geq 0$}. By induction it follows that
$$\Omega^{np^k} - \Omega^{p^k} \in p^{k+1}A(S,S)$$
for \mbox{$k \geq 0,$} \mbox{$n > 0$}. In particular
$$\Omega^{p^l} - \Omega^{p^k} \in p^{k+1}A(S,S)$$
when \mbox{$l \geq k > 0$}, so
$$\Omega, \Omega^p, \Omega^{p^2}, \dots$$
is a Cauchy sequence in the $p$--adic topology of $A(S,S)$. Hence
it converges to a unique element \mbox{$\CharIdem \in
\pComp{A(S,S)}$}. Since the multiplication in $A(S,S)$ is
continuous with respect to the $p$--adic topology, ${\CharIdem}^2$
is the limit of the sequence
$$\Omega^2, \Omega^{2p}, \Omega^{2p^2}, \dots$$
Idempotence of ${\CharIdem}$ now follows by taking the limit of
\ref{eq:Convergent1} over $k.$
\end{proof}

We can now prove the main result of this section.
\begin{proposition} \label{prop:newIdempotent}
Every saturated fusion system has a characteristic idempotent.
\end{proposition}
\begin{proof}
Let $\F$ be a saturated fusion system over a finite $p$--group $S$.
Take a characteristic $(S,S)$--biset $\Omega$ as given by 
\fullref{prop:BLOIdempotent}. By Lemmas \ref{lem:ExistsM} and
\ref{lem:CanTakePowers} we may assume that $\Omega$ is idempotent
$mod~ p$. By \fullref{lem:Convergent} the sequence
$$\Omega, \Omega^p, \Omega^{p^2}, \dots$$
converges to an idempotent \mbox{$\CharIdem \in \pComp{A(S,S)}$}. We
show that $\CharIdem$ has the Linckelmann--Webb properties.

By an induction similar to that in 
\fullref{lem:IncreasingIdempotents} one can show that
\mbox{$\countS(\Omega) \equiv 1~ (mod~ p)$} implies that
\mbox{$\countS(\Omega^{p^k}) \equiv 1~ (mod~ p^{k+1})$} for
\mbox{$k \geq 0$}. It follows that \mbox{$\countS(\CharIdem) =
1$}, proving (c).

It is not hard to see that $A_\F(S,S)$ is a closed subspace of
$A(S,S)$ in the $p$--adic topology and hence that
$\pComp{A_{\F}(S,S)}$ is a closed subspace of $\pComp{A(S,S)}$.
Since each $\Omega^n$ is in $A_{\F}(S,S)$ by 
\fullref{cor:RFclosed}, it follows that the limit ${\CharIdem}$ is in
$\pComp{A_{\F}(S,S)}$, proving (a).

Let \mbox{$P\leq S$} and \mbox{$\varphi \in \HomF{P}{S}$}. By
Property $\rm(b1')$ we have
$$\Omega \circ [P,\varphi]_P^S = \Omega \circ [P,\iota_P]_P^S$$
and consequently
$$\Omega^{p^k} \circ [P,\varphi]_P^S = \Omega^{p^k} \circ [P,\iota_P]_P^S,$$
for all \mbox{$k \geq 0$}. Since the pairing
$$\circ\negmedspace: A(S,S) \times A(P,S) \to A(P,S)$$ is
continuous in the $p$--adic topology on the relevant
$\mathbb{Z}$--modules, we can take limits to get
$${\CharIdem} \circ [P,\varphi]_P^S = {\CharIdem} \circ [P,\iota_P]_P^S,$$
proving (b1). Property (b2) follows similarly from Property $\rm(b2')$.
\end{proof}

\section{Properties of characteristic idempotents} \label{sec:Inherent}
In this section we perform a further study of the characteristic
idempotents introduced in the previous section. We discover that
the effect of multiplicaton by a characteristic idempotent on
$\pComp{A(S_1,S_2)}$ essentially amounts to quotienting out fusion
conjugacy in the source or target, as appropriate. This allows us
to glean important information about the structure of a
characteristic idempotent in the $[P,\varphi]$--basis, which allows
us to prove its uniqueness, and will also prove surprisingly
useful for proving later naturality results.

\begin{proposition} \label{prop:omegaOnConjugates}
Let $\F_1$ and $\F_2$ be saturated fusion systems over the finite
\mbox{$p$--groups} $S_1$ and $S_2$, respectively, and let
$\CharIdem_1$ and $\CharIdem_2$ be characteristic idempotents of
$\F_1$ and $\F_2$, respectively. If the $(S_1,S_2)$--pairs $(P,\psi)$
and $(Q,\rho)$ are $(\F_1,\F_2)$--conjugate, then
$$\CharIdem_2 \circ [Q,\rho] \circ \CharIdem_1 = \CharIdem_2 \circ [P,\psi] \circ \CharIdem_1.$$
\end{proposition}
\begin{proof}
Since $\CharIdem_1$ is left $\F_1$--stable and $\CharIdem_2$ is
right $\F_2$--stable, this is a special case of 
\fullref{lem:FinvOnConjugates}.
\end{proof}

Although the following proposition may be of limited interest in
its own right, it is the central result of this paper.
\begin{proposition} \label{prop:omegaBasis}
Let $\F_1$ and $\F_2$ be saturated fusion systems over the finite
\mbox{$p$--groups} $S_1$ and $S_2$, respectively and let
$\CharIdem_1$ and $\CharIdem_2$ be characteristic idempotents of
$\F_1$ and $\F_2$, respectively. Let $I$ be the set of
$(\F_1,\F_2)$--conjugacy classes of $(S_1,S_2)$--pairs, and pick a
representative $(P_i,\psi_i)$ for each \mbox{$i \in I$}. Then the
collection
$$\{\CharIdem_2 \circ [P_i,\psi_i] \circ \CharIdem_1 \mid i\in I \}  $$
forms a $\Zp$--basis for \mbox{$\CharIdem_2 \circ
\pComp{A(S_1,S_2)} \circ \CharIdem_1$}.
\end{proposition}
\begin{proof}
It follows from \fullref{prop:omegaOnConjugates} that
\mbox{$\CharIdem_2 \circ \pComp{A(S_1,S_2)} \circ \CharIdem_1$} is
spanned by the collection, so it suffices to prove linear
independence. By Property (b) of characteristic idempotents and
\fullref{lem:PreservesFiltration} we have
$$\CharIdem_2 \circ [P_i,\psi_i] \circ \CharIdem_1 \in \noFsubModule{P_i}{\psi_i}$$
for each \mbox{$i \in I$}. Note however, that
$$\countS\left(\CharIdem_2 \circ [P_i,\psi_i] \circ \CharIdem_1\right) = \countS(\CharIdem_2) \cdot \countS([P_i,\psi_i]) \cdot \countS(\CharIdem_1) = 1 \cdot |S_1/P_i| \cdot 1 = |S_1/P_i|,$$
whereas
\begin{equation} \label{eq:strictsubAugmentation}
 \countS\left(\noFstrictsubModule{P_i}{\psi_i}\right) \subseteq p~|S_1/P_i|~\Zp.
\end{equation}
Therefore,
\begin{equation} \label{eq:omegaPreservesLayer}
  \CharIdem_2 \circ [P_i,\psi_i] \circ \CharIdem_1 \in \noFsubModule{P_i}{\psi_i} \setminus \noFstrictsubModule{P_i}{\psi_i}
\end{equation}
 for each
\mbox{$i \in I$}.

Now, let \mbox{$c_i \in \Zp$} for each \mbox{$i \in I$} and assume
that
\begin{equation} \label{eq:LinearDependence}
  \sum_{i \in I} c_i \cdot \left(\CharIdem_2 \circ [P_i,\psi_i] \circ \CharIdem_1\right) = 0.
\end{equation}
Put
$$I' = \{i \in I \mid c_i \neq 0\}.$$
If $I'$ is nonempty, then let $j$ be a maximal element of $I'$
regarded as a poset under $(\F_1,\F_2)$--subconjugacy. By
\ref{eq:omegaPreservesLayer} there is a $(S_1,S_2)$--pair
\mbox{$(Q,\rho) \noFcon (P_j,\psi_j)$} such that
$$\chi_{[Q,\rho]}\left(  \CharIdem_2 \circ [P_j,\psi_j] \circ \CharIdem_1   \right) \neq 0.$$
On the other hand, for \mbox{$i \in I'\setminus \{j\}$} the
maximality of $j$ implies that $[Q,\rho]$ is not
$(\F_1,\F_2)$--subconjugate to $(P_i,\psi_i)$. Hence
$$\chi_{[Q,\rho]}\left( \noFsubModule{P_i}{\psi_i} \right) = 0$$
and in particular
$$\chi_{[Q,\rho]}\left(  \CharIdem_2 \circ [P_i,\psi_i] \circ \CharIdem_1   \right) = 0.$$
Now we get
\begin{align*}
  \chi_{[Q,\rho]}\left(\sum_{i \in I} c_i \cdot \left(\CharIdem_2 \circ [P_i,\psi_i] \circ \CharIdem_1 \right)\right)
  &= \sum_{i \in I} c_i \cdot \chi_{[Q,\rho]}\left( \CharIdem_2 \circ [P_i,\psi_i] \circ \CharIdem_1 \right)\\
  &= \sum_{i \in I\setminus I'} \underbrace{c_i}_{=0} \cdot \chi_{[Q,\rho]}\left( \CharIdem_2 \circ [P_i,\psi_i] \circ \CharIdem_1 \right)\\
  &+ \sum_{i \in I'\setminus\{j\}} c_i \cdot \underbrace{\chi_{[Q,\rho]}\left( \CharIdem_2 \circ [P_i,\psi_i] \circ \CharIdem_1 \right)}_{=0}\\
  &+ \underbrace{c_j}_{\neq 0} \cdot \underbrace{\chi_{[Q,\rho]}\left( \CharIdem_2 \circ [P_j,\psi_j] \circ \CharIdem_1 \right)}_{\neq 0}\\
  &\neq 0,
\end{align*}
contradicting \ref{eq:LinearDependence}. Therefore $I'$ must be
empty and we conclude that the collection is linearly independent.
\end{proof}

\begin{remark} \label{rem:SingleomegaBasis}
Since the multiplicative identity $[S_1,id]$ of $\pComp{A(S_1,S_1)}$
is a characteristic idempotent for $\F_{S_1}$ and $[S_2,id]$ is a
characteristic idempotent for $\F_{S_2}$, Propositions
\ref{prop:omegaOnConjugates} and \ref{prop:omegaBasis} can also be
used to obtain a basis for \mbox{$\CharIdem_2 \circ
\pComp{A(S_1,S_2)}$} and \mbox{$\pComp{A(S_1,S_2)} \circ
\CharIdem_1$}.
\end{remark}

We determine the structure of a characteristic idempotent
$\CharIdem$ by carefully analyzing the idempotence relation
\mbox{$\CharIdem \circ \CharIdem = \CharIdem$}. In light of the
previous proposition, a convenient tool for doing this is the
projection
$$\StableFusion{\F}{S} \longrightarrow \CharIdem \StableFusion{\F}{S},$$
given by left multiplication by $\CharIdem$. This projection can be
easily described by
$$\Omega = \sum_{[P,\varphi] \noFsub [S,id]} \chi_{[P,\varphi]}(\Omega) \cdot [P,\varphi] \longmapsto  \CharIdem \circ \Omega = \sum_{[P] \noFsub [S]} \chi_{[P]}^{\F}(\Omega) \cdot (\CharIdem \circ [P,\iota_P]),$$
where the homomorphisms $\chi_{[P]}^{\F}$ are as in the following
definition.
\begin{definition}
Let $\F$ be a fusion system over a finite $p$--group $S$. For each
$S$--conjugacy class $[P]$ of subgroups of $S$, let $\chi_{[P]}^{\F}$
and $\chi^{[P]}_{\F}$ be the homomorphisms \mbox{$\pComp{A(S,S)} \to
\Zp$} given by
$$\chi_{[P]}^{\F} = \sum_{[P,\varphi] \Ffcon{(\F_S,\F)} [P,\iota_P]} \chi_{[P,\varphi]} $$
and
$$\chi^{[P]}_{\F} = \sum_{[Q,\varphi] \Ffcon{(\F,\F_S)} [P,\iota_P]} \chi_{[Q,\varphi]}.$$
\end{definition}
Note that similarly we have
$$\Omega \circ \CharIdem = \sum_{[P] \noFsub [S]} \chi^{[P]}_{\F}(\Omega) \cdot ([P, \iota_P] \circ \CharIdem). $$
The following lemma now effectively allows us to determine the
structure of characteristic idempotents.
\begin{lemma}\label{lem:chiSumsToZero}
Let $\F$ be a saturated fusion system over a finite $p$--group $S$
and let $\CharIdem$ be a characteristic idempotent of $\F$. Then
$$\CharIdem = \CharIdem^{\noFfcon} + \CharIdem^{\noFstrictsub},$$
where
$$\CharIdem^{\noFfcon} = \frac{1}{|\Out{S}{\F}|}\sum_{\varphi \in \Out{S}{\F}}[S,\varphi] \in \FconModule{\F}{\F}{S}{id}$$
and
$$\CharIdem^{\noFstrictsub} \in \FstrictsubModule{\F}{\F}{S}{id}.$$
Furthermore,
$$\chi_{[P]}^{\F}(\CharIdem) = 0$$
and
$$\chi^{[P]}_{\F}(\CharIdem) = 0$$
for all proper subgroups \mbox{$P < S$}.
\end{lemma}
\begin{proof}
There is a direct sum of modules
$$\StableFusion{\F}{S} = \noFsubModule{S}{id} = \noFconModule{S}{id} \oplus \noFstrictsubModule{S}{id},$$
where $\noFstrictsubModule{S}{id}$ is a two-sided ideal of
$\StableFusion{\F}{S}$ (by \fullref{lem:PreservesFiltration}), and
$\noFconModule{S}{id}$ is a subring of $\StableFusion{\F}{S}$ (as
can be easily checked). We can therefore uniquely write
$$\CharIdem = \CharIdem^{\noFfcon} + \CharIdem^{\noFstrictsub},$$
where
$$\CharIdem^{\noFfcon} \in \noFconModule{S}{id}$$
and
$$\CharIdem^{\noFstrictsub} \in \noFstrictsubModule{S}{id}.$$
From Properties (a) and (b) of characteristic idempotents, one
deduces by standard techniques that $\CharIdem^{\noFfcon}$ must be
of the form
$$\CharIdem^{\noFfcon} = \alpha\sum_{\varphi \in \Out{S}{\F}}[S,\varphi],$$
where \mbox{$\alpha \in \Zp$}. Now,
$$\CharIdem^{\noFfcon} + \CharIdem^{\noFstrictsub} = \CharIdem = \CharIdem^2 = (\CharIdem^{\noFfcon})^2 + (\CharIdem^{\noFfcon} \circ \CharIdem^{\noFstrictsub} + \CharIdem^{\noFstrictsub} \circ \CharIdem^{\noFfcon} + (\CharIdem^{\noFstrictsub})^2),$$
where \mbox{$(\CharIdem^{\noFfcon})^2 \in \noFconModule{S}{id}$}
and the second term is in the ideal $\noFstrictsubModule{S}{id}$.
By uniqueness of such decompositions $\CharIdem^{\noFfcon}$ is
therefore an idempotent, and since
$$(\CharIdem^{\noFfcon})^2 = \alpha \cdot |\Out{S}{\F}| \cdot \CharIdem^{\noFfcon},$$
this implies that
$$\alpha \in \{0,|\Out{S}{\F}|^{-1}\}.$$
(Note that $|\Out{S}{\F}|$ is prime to $p$ since $\F$ is
saturated.) By Property (c) of characteristic idempotents we have
$$1 = \countS(\CharIdem) = \countS(\CharIdem^{\noFfcon}) + \countS(\CharIdem^{\noFstrictsub}).$$
As a special case of \ref{eq:strictsubAugmentation} we have
\mbox{$\countS(\noFstrictsubModule{S}{id}) \subseteq p~\Zp$}, and it
follows that
$$\alpha \cdot |\Out{S}{\F}| = \countS(\CharIdem^{\noFfcon}) \equiv \countS(\CharIdem) \equiv 1~ (mod~ p).$$
Hence \mbox{$\alpha \neq 0$}, leaving only the possibility
$$\CharIdem^{\noFfcon} = \frac{1}{|\Out{S}{\F}|}\sum_{\varphi \in \Out{S}{\F}}[S,\varphi].$$
To prove the second claim we start by deducing from Property (b)
and the description of $\CharIdem^{\noFfcon}$ above, that
$$\CharIdem \circ \CharIdem^{\noFfcon} = \CharIdem.$$
Hence the idempotence of $\CharIdem$ yields
$$\CharIdem = \CharIdem \circ \CharIdem = \CharIdem \circ (\CharIdem^{\noFfcon} + \CharIdem^{\noFstrictsub}) = \CharIdem + \CharIdem \circ \CharIdem^{\noFstrictsub},$$
and we get
$$\CharIdem \circ \CharIdem^{\noFstrictsub} = 0.$$
By Property (a) and \fullref{rem:RFDescription}, we can write
$$\CharIdem^{\noFstrictsub} = \sum_{[P,\varphi] \noFstrictsub [S,id]} \chi_{[P,\varphi]}(\CharIdem) \cdot [P,\varphi],$$
so
$$ 0 = \CharIdem \circ \CharIdem^{\noFstrictsub} = \sum_{[P] \noFstrictsub [S]} \chi_{[P]}^{\F}(\CharIdem) \cdot \left(\CharIdem\circ[P,\iota_P]\right),$$
and the result for $\chi_{[P]}^{\F}(\CharIdem)$ follows upon
noting that the collection
$$\{\CharIdem \circ [P,\iota_P] \mid [P] \noFstrictsub [S] \}$$
is linearly independent over $\Zp$ by 
\fullref{prop:omegaBasis} and \fullref{rem:SingleomegaBasis}.

The result for $\chi^{[P]}_{\F}(\CharIdem)$ is proved similarly by
performing the analogous simplifications of \mbox{$
\CharIdem^{\noFstrictsub} \circ \CharIdem= 0$}. We omit the details.
\end{proof}

The lemma has an interesting consequence.
\begin{proposition} \label{prop:omegaUnique}
Every saturated fusion system has a unique characteristic
idempotent.
\end{proposition}
\begin{proof}
Let $\F$ be a saturated fusion system over a finite $p$--group $S$.
By \fullref{prop:newIdempotent}, $\F$ has a characteristic
idempotent $\CharIdem$. We proceed to prove uniqueness. Recalling
that we can write
$$\omega = \sum_{[P,\varphi] \noFsub [S,id] } \chi_{[P,\varphi]}(\CharIdem) \cdot [P,\varphi],$$
the method of proof is to set up a system of linear equations for
the coefficients $\chi_{[P,\varphi]}(\CharIdem)$, and show that this
system is fully determined and thus uniquely determines the
coefficients. Note that we need not show that the system is solvable
since the existence of a characteristic idempotent has already been
established.

To present the argument it is helpful to index the basis of
$\StableFusion{\F}{S}$. Take representatives \mbox{$P_0, P_1,
\ldots, P_n$} for the $S$--conjugacy classes of subgroups of $S$,
labelled in such a way that
$$|P_0| \geq |P_1| \geq \ldots \geq |P_n|.$$
For each \mbox{$i \in \{0,1,\ldots,n\},$} pick representatives
\mbox{$\varphi_{i0}, \varphi_{i1}, \ldots, \varphi_{im_i}$} for the
$S$--conjugacy classes of homomorphisms in $\HomF{P}{S}$, now
labelled so that \mbox{$\varphi_{i0} = \iota_{P_i}$}. The collection
$\{ [P_i,\varphi_{ij}] \mid 0 \leq i \leq n, 0 \leq j \leq m_i \}$
is then a basis for $\StableFusion{\F}{S}$. We order the basis
through the lexicographic order of the indexing set:
$$(0,0) \leq (0,1) \leq \cdots \leq (0,m_0) \leq (1,0) \leq (1,1) \leq \cdots \leq (n,m_n).$$
This ordering has the property that \mbox{$[P_i,\varphi_{ij}]
\Fsub{(S,S)} [P_k,\varphi_{kl}]$} implies \mbox{$(k,l) \leq (i,j).$}

For the remainder of this proof, a pair $(i,j)$ is understood to
satisfy \mbox{$0 \leq i \leq n$} and \mbox{$1 \leq j \leq m_i$}. For
such a pair $(i,j)$, write
$$\Graph{i}{j} := \Graph{P_i}{\varphi_{ij}}$$
and
$$c_{ij} :=  \chi_{[P_i,\varphi_i^j]} (\CharIdem)$$
for short. We will also use the shorthand notation $\sum_{i,j}$
for the double sum $\sum_{i=0}^n \sum_{j=0}^{m_i}.$ Using this
notation we can write
$$\CharIdem = \sum_{i,j} c_{ij} [P_i,\varphi_{ij}].$$
Note that by \fullref{lem:StableFxdPts}, Property (b1) implies
that
\begin{equation*} %\label{eq:FxdPts}
 |\CharIdem^{\Graph{i}{j}}| = |\CharIdem^{\Graph{i}{0}}|
\end{equation*}
for all pairs $(i,j)$. This can be rewritten as the equation
\begin{equation*} %\label{eq:ij}
  \sum_{k,l} c_{kl} \cdot (|[P_k,\varphi_{kl}]^{\Graph{i}{j}}| - |[P_k,\varphi_{kl}]^{\Graph{i}{0}}|) = 0,
\end{equation*}
which we refer to as Equation $(i,j)$ when \mbox{$j \neq 0$}. When
\mbox{$j = 0$}, this equation becomes trivial. Instead we consider
the equations given by \fullref{lem:chiSumsToZero}. That is, we
let Equation $(0,0)$ be the equation
$$\sum_{j = 0}^{m_0} c_{0j} = 1,$$
and for \mbox{$i = 1, \ldots, n$}, we let Equation $(i,0)$ be the
equation
$$\sum_{j = 0}^{m_i} c_{ij} = 0.$$
If we now write Equations $(0,0)$ to $(i,m_i)$ one below the other
going by the lexicographic order, we obtain a system of equations,
which can be represented on matrix form as
$$A c = b,$$
where $c$ is a vector with entries $c_{ij}$, $b$ is a vector with
$1$ in its first entry and $0$ everywhere else, while $A$ is a
square matrix whose rows and columns are both indexed by pairs
$(i,j)$. The proof is completed by showing that $A$ has nonzero
determinant.

There is an obvious way to regard $A$ as a $(n \times n)$ block
matrix where the blocks are indexed by $i$. We show that $A$ is in
fact a lower triangular block matrix. Since
$$|[Q,\psi] ^ {\Graph{R}{\rho}}| = | \Graph{Q}{\psi} \backslash  N_{S\times S}( \Graph{R}{\rho},\Graph{Q}{\psi} ) | = \frac{|N_{S\times S}( \Graph{R}{\rho},\Graph{Q}{\psi} ) |}{ | Q |}$$
for $(S,S)$--pairs $(Q,\psi)$ and $(R,\rho)$, we see by 
\fullref{rem:SubconjugacyLink} that
$$|[Q,\psi] ^ {\Graph{R}{\rho}}| \neq 0 $$
if and only if
$$[R,\rho] \Fsub{(S,S)} [Q,\psi].$$
The chosen order of the basis of $\StableFusion{\F}{S}$ therefore
guarantees that $A$ is a lower triangular block matrix.

It now suffices to show that the matrices occurring on the
diagonal of $A$ have nonzero determinant. The $i$-th matrix on the
diagonal has the form
$$A_i = \left(
          \begin{array}{ccccccc}
            1       & 1      & 1      &        & 1          & 1          & 1 \\
            -a_{i0} & a_{i1} & 0      & \dots  & 0          & 0          & 0 \\
            -a_{i0} & 0      & a_{i2} &        & 0          & 0          & 0 \\
                    & \vdots &        & \ddots &            & \vdots     &   \\
            -a_{i0} & 0      & 0      &        & a_{im_i-2} & 0          & 0 \\
            -a_{i0} & 0      & 0      & \dots  & 0          & a_{im_i-1} & 0 \\
            -a_{i0} & 0      & 0      &        & 0          & 0          & a_{im_i} \\
          \end{array}
        \right)
$$
where
$$a_{ij} = \frac{|N_{S \times S}(\Graph{i}{j})|}{|P_i|}.$$
Direct calculation shows that the determinant of $A_i$ is
$$\det(A_i) = \sum_{j=0}^{m_i} \frac{\prod_{k=0}^{m_i}a_{ik}}{a_{ij}},$$
which is nonzero since all the $a_{ij}$'s are positive integers.
\end{proof}

The previous proposition allows us to speak of \emph{the}
characteristic idempotent of a saturated fusion system.
\begin{definition}
For a saturated fusion system $\F$, let $\CharIdem_{\F}$ denote the
characteristic idempotent of $\F$.
\end{definition}

\begin{remark} We make the following observations about the proof of 
\fullref{prop:omegaUnique}.
\begin{enumerate}
  \item
As a byproduct of the proof we have produced an explicit system of
equations which we can solve to produce $\CharIdem_{\F}$.
  \item
The coefficients in these equations are all integers. Therefore
$\CharIdem_{\F}$ can be regarded as an element in the
$p$--localized double Burnside ring:
$$\CharIdem_{\F} \in A(S,S)_{(p)} = A(S,S) \otimes \Z_{(p)}.$$
  \item
The proof actually shows that $\CharIdem_{\F}$ is the unique
idempotent in $\pComp{A(S,S)}$ (or $A(S,S)_{(p)}$) with the
Linckelmann--Webb Properties (a), (b1) and (c). A symmetric
argument shows that $\CharIdem_{\F}$ is the unique idempotent in
$\pComp{A(S,S)}$ with Properties (a), (b2) and (c). We are
therefore in the surprising situation that for an idempotent in
$\pComp{A(S,S)}$ with Properties (a) and (c), the presence of
either stability Property (b1) or (b2) implies the presence of the
other.
\end{enumerate}
\end{remark}

\section{Stable characteristic idempotents} \label{sec:StableCharIdempotent}
In the preceding sections we have undertaken a thorough study of
the structure and properties of characteristic idempotents of
saturated fusion systems. Going by \cite{BLO2}, we will later
construct the classifying spectrum of a saturated fusion system
$\F$ over a finite $p$--group $S$ as the mapping telescope of the
stable characteristic idempotent $\StableCharIdem_\F$ of
$\ClSpectrum{S}$ induced by the characteristic idempotent
$\CharIdem_\F$. In this section we briefly study the stable
characteristic idempotent and interpret the consequences of the
results of the preceding two sections in terms of stable maps,
with a particular view towards fusion stability properties as
defined below.

\begin{definition}
Let $\F$ be a fusion system over a finite group $S$ and let $X$ be
a spectrum. We say a map \mbox{$f\negmedspace: \ClSpectrum{S} \to
X$} is \emph{right $\F$--stable} if for every \mbox{$P \leq S$} and
\mbox{$\varphi \in \HomF{P}{S}$} we have
$$f\circ \ClSpectrum{\varphi} \simeq f \circ \ClSpectrum{\iota_P}.$$
Similarly we say that a map \mbox{$f\negmedspace: X \to
\ClSpectrum{S}$} is \emph{left $\F$--stable} if for every \mbox{$P
\leq S$} and \mbox{$\varphi \in \HomF{P}{S}$} we have
$$tr_{\varphi} \circ f \simeq tr_P \circ f.$$
\end{definition}
We will often simply say that a map is $\F$--stable. Unless
\mbox{$X \simeq \ClSpectrum{S}$} in the previous definition, it
should be clear from the context whether this means left or right
$\F$--stability.

In \fullref{sec:Burnside} we defined
$\pComp{\wtilde{A}(S,S)}$ to be the quotient module of
$\pComp{A(S,S)}$ obtained by quotienting out all trivial basis
elements $[P,\psi]$, and discussed the isomorphism
$$\wtilde{\alpha} : \pComp{\wtilde{A}(S,S)} \stackrel{\cong}{\longrightarrow} \StableMaps{BS}{BS} $$
of the Segal Conjecture. Since $\pComp{\wtilde{A}(S,S)}$ is a
free $\Zp$--module with one basis element $[P,\psi]$ for each
conjugacy class of nontrivial $(S,S)$--pairs, we can also regard
$\pComp{\wtilde{A}(S,S)}$ as a submodule (but not subring) of
$\pComp{A(S,S)}$ and we do so in this section. From this point of
view we define a version of the augmentation $\countS$ by
$$\wtilde{\countS} := \countS \circ \wtilde{\alpha}^{-1} \negmedspace : \StableMaps{BS}{BS} \stackrel{\wtilde{\alpha}^{-1}}{\longrightarrow} \pComp{\wtilde{A}(S,S)} \subseteq \pComp{A(S,S)} \stackrel{\countS}{\longrightarrow} \Zp.$$
This is not a ring homomorphism and therefore not an augmentation,
although it plays a similar role.
\begin{definition}
Let $\F$ be a fusion system over a finite $p$--group $S$. A
\emph{stable characteristic idempotent for $\F$} is an idempotent
\mbox{$\CharIdem \in \StableMaps{BS}{BS}$} with the following
properties:
\begin{enumerate}
  \item[(a)]
$\StableCharIdem \in \wtilde{\alpha}(\StableFusion{\F}{S})$,
  \item[(b1)]
$\StableCharIdem$ is right $\F$--stable.
  \item[(b2)]
$\StableCharIdem$ is left $\F$--stable.
  \item[(c)]
$\wtilde{\countS}(\StableCharIdem) = 1$.
\end{enumerate}
\end{definition}

\begin{proposition}
Every saturated fusion system $\F$ has a unique stable
characteristic idempotent $\StableCharIdem_{\F}$.
\end{proposition}
\begin{proof}
If $\F$ is a saturated fusion system over a finite $p$--group $S$,
then the image $\StableCharIdem_{\F}$ in $\StableMaps{BS}{BS}$ of
the characteristic idempotent $\CharIdem_{\F}$ under the
surjection
\[
\xymatrix{
  \pComp{A(S,S)} \ar@{->>}[0,1]  & \pComp{\wtilde{A}(S,S)}  \ar[0,1]^{\wtilde{\alpha}}_{\cong} &  \StableMaps{BS}{BS}
}
\]
is a stable characteristic idempotent for $\F$.

To prove uniqueness, suppose that $\StableCharIdem$ is a stable
characteristic idempotent of $\F$ and consider \mbox{$\CharIdem :=
\wtilde{\alpha}^{-1}(\StableCharIdem)$}, regarded as an element
of $\pComp{A(S,S)}$. It is enough to show that $\omega$ is a
characteristic idempotent for $\F$, because the result then
follows from uniqueness of characteristic idempotents (\fullref{prop:omegaUnique}).

It is immediate that \mbox{$\CharIdem \in \StableFusion{\F}{S}$} and
\mbox{$\countS(\CharIdem) = 1$}. We show that $\CharIdem$ is
idempotent. Indeed, idempotence of $\StableCharIdem$ implies that
$$\nu := \CharIdem \circ \CharIdem - \CharIdem$$
is a linear combination of trivial basis elements. But since $\nu$
is an element of $\StableFusion{\F}{S}$, and the only trivial basis
element in $\StableFusion{\F}{S}$ is $[1,triv]$, where $1$ denotes
the trivial subgroup, $\nu$ must be of the form \mbox{$\nu = a \cdot
[1,triv]$}. Applying $\countS$ we get
$$1 = \countS(\omega \circ \omega) = \countS(\omega + \nu) = 1 + a \cdot |S|.$$
We conclude that \mbox{$a = 0$}, whence \mbox{$\nu = 0$} and
$\omega$ is idempotent.

It only remains to prove that $\CharIdem$ has Property (b). This
is deduced from Property (b) for $\StableCharIdem$ in a similar
way.
\end{proof}

As a consequence of \fullref{lem:chiSumsToZero} we now show that
stable characteristic idempotents completely characterize
fusion-stable maps.
\begin{corollary} \label{cor:FinvEatsOmega}
Let $\F$ be a saturated fusion system over a finite $p$--group $S$
and let $X$ be a spectrum. A map \mbox{$f: \ClSpectrum{S} \to X $}
is $\F$--stable if and only if
$$f \circ \StableCharIdem_{\F} \simeq f.$$
Similarly, a map \mbox{$f\negmedspace : X \to \ClSpectrum{S}$} is
$\F$--stable if and only if
$$\StableCharIdem_{\F} \circ f \simeq f.$$
\end{corollary}
\begin{proof}
Suppose first that \mbox{$f: \ClSpectrum{S} \to X $} is
$\F$--stable. Note that we can write
$$\CharIdem_{\F}  \simeq  \sum_{[P] \noFsub [S]} \left( \sum_{[P,\varphi] \Fcon{(\F_{S},\F)} [P,\iota_P]} \chi_{[P,\varphi]}(\CharIdem_{\F}) [P,\varphi] \right) $$
Now,
\begin{align*}
  f \circ \StableCharIdem_{\F}
  &{\simeq} f \circ \sum_{[P] \noFsub [S]} \left( \sum_{[P,\varphi] \Fcon{(\F_{S},\F)} [P,\iota_P]} \chi_{[P,\varphi]}(\CharIdem_{\F}) \cdot \wtilde{\alpha}([P,\varphi]) \right)\\
  &{\simeq} f \circ \sum_{[P] \noFsub [S]} \left( \sum_{[P,\varphi] \Fcon{(\F_{S},\F)} [P,\iota_P]} \chi_{[P,\varphi]}(\CharIdem_{\F}) \cdot \left(\ClSpectrum{\varphi} \circ tr_P\right) \right)\\
  &{\simeq} \sum_{[P] \noFsub [S]} \left( \sum_{[P,\varphi] \Fcon{(\F_{S},\F)} [P,\iota_P]} \chi_{[P,\varphi]}(\CharIdem_{\F}) \cdot \left(f \circ \ClSpectrum{\varphi} \circ tr_P \right) \right)\\
  &{\simeq} \sum_{[P] \noFsub [S]} \left( \left( \sum_{[P,\varphi] \Fcon{(\F_{S},\F)} [P,\iota_P]} \chi_{[P,\varphi]}(\CharIdem_{\F}) \right) \cdot \left(f \circ \ClSpectrum{\iota_P} \circ tr_P \right) \right)\\
  &{\simeq} \sum_{[P] \noFsub [S]} \left( \chi_{[P]}^{\F}(\CharIdem_{\F}) \cdot \left(f \circ \ClSpectrum{\iota_P} \circ tr_P \right) \right).
\end{align*}
By \fullref{lem:chiSumsToZero}, we have
$$\chi_{[P]}^{\F}(\CharIdem_{\F}) = \begin{cases}
                                   {1} & \text{if $P = S$},\\
                                   {0} & \text{if $P < S$}.
                                \end{cases}
$$
We conclude that
$$f \circ \StableCharIdem_{\F} \simeq f \circ \ClSpectrum{\iota_S} \circ tr_S  \simeq  f.$$
The converse implication follows from $\F$--stability of
$\StableCharIdem_{\F}$.

The result for left $\F$--stable maps is deduced similarly by
writing
$$\CharIdem_{\F} = \sum_{[P] \noFsub [S]} \left( \sum_{[Q,\varphi] \Fcon{(\F,\F_{S})} [P,\iota_P]} \chi_{[Q,\varphi]}(\CharIdem_{\F}) [Q,\varphi] \right).$$
If $f$ is a left $\F$--stable map $f$,  we get the following after
some manipulations
$$\StableCharIdem_{\F} \circ f \simeq \sum_{[P] \noFsub [S]} \left( \chi^{[P]}_{\F}(\CharIdem_{\F}) \cdot \left(\ClSpectrum{\iota_P} \circ tr_P \circ f\right) \right)$$
and \mbox{$\StableCharIdem_{\F} \circ f \simeq f$} follows after
appealing to \fullref{lem:chiSumsToZero}.
\end{proof}

\section{Classifying spectra} \label{sec:ClSpectrum}
In this section we develop the theory of classifying spectra for
saturated fusion systems. The classifying spectra constructed here
agree with those suggested in \cite{BLO2}, but the added rigidity of
the new construction allows us to prove that the assignment is
functorial for fusion-preserving homomorphisms between fusion
systems. We will also endow the classifying spectrum of a saturated
fusion system $\F$ over a finite $p$--group $S$ with additional
structure by regarding it as an object under $\ClSpectrum{S}$. We
support this idea by proving that we can reconstruct the fusion
system $\F$ from its classifying spectrum, when regarded as an
object under $\ClSpectrum{S}$.
\smallskip

Let $\F$ be a saturated fusion system over a finite $p$--group $S$
and denote the infinite mapping telescope of the stable
characteristic idempotent $\StableCharIdem_{\F}$ by
$\ClSpectrum{\F}$. In other words,
$$\ClSpectrum{\F} := \hocolim{} \left(\ClSpectrum{S} \stackrel{\StableCharIdem_{\F}}{\longrightarrow} \ClSpectrum{S} \stackrel{\StableCharIdem_{\F}}{\longrightarrow} \ClSpectrum{S} \stackrel{\StableCharIdem_{\F}}{\longrightarrow} \dots \right).$$
We denote the structure map of the homotopy colimit by
\mbox{$\sigma_{\F}\negmedspace: \ClSpectrum{S} \to
\ClSpectrum{\F}$}. Since $\StableCharIdem_{\F}$ is idempotent up
to homotopy, we get a homotopy factorization of
$\StableCharIdem_{\F}$ through the homotopy colimit
$$\StableCharIdem_{\F}\negmedspace: \ClSpectrum{S} \stackrel{\sigma_{\F}}{\longrightarrow} \ClSpectrum{\F} \stackrel{t_{\F}}{\longrightarrow} \ClSpectrum{S},$$
such that \mbox{$\sigma_{\F} \circ t_{\F} \simeq
id_{\ClSpectrum{\F}}$}. Note that $t_\F$ is, up to homotopy, the
unique map with these properties. Since $\ClSpectrum{\F}$ is a
retract of the $p$--complete spectrum $\ClSpectrum{S}$, it is
\mbox{$p$--complete}.

\begin{definition}
Let $\F$ be a saturated fusion system over a finite $p$--group $S$.
The \emph{classifying spectrum of $\F$} is the spectrum
$\ClSpectrum{\F}$, the \emph{structure map of $\F$} is the map
$\sigma_{\F}$, and the \emph{transfer of $\F$} is the map
$$t_{\F}\negmedspace: \ClSpectrum{\F} \to \ClSpectrum{S},$$
such that \mbox{$\sigma_{\F} \circ t_{\F} \simeq
id_{\ClSpectrum{\F}}$} and \mbox{$t_{\F} \circ \sigma_{\F} \simeq
\StableCharIdem_{\F}$}.

The \emph{structured classifying spectrum of $\F$} is the mapping
telescope
$$\sigma_{\F} \negmedspace: \ClSpectrum{S} \to \ClSpectrum{\F},$$
regarded as an object under $\ClSpectrum{S}$.
\end{definition}
In the special case when \mbox{$\F = \F_S$} is the fusion system of
$S$, the stable characteristic idempotent $\StableCharIdem_{\F_S}$
is just the identity of $\ClSpectrum{S}$, so the structured
classifying spectrum of $\F_S$ becomes a natural homotopy
equivalence
$$\sigma_{\F_S} : \ClSpectrum{S} \stackrel{\simeq}{\longrightarrow} \ClSpectrum{\F_S}.$$
We will therefore often replace $\ClSpectrum{\F_S}$ by
$\ClSpectrum{S}$.

It has been shown in \cite{BLO2} that for a $p$--local finite group
$\plfg$, the classifying spectrum $\ClSpectrum{\F}$ is homotopy
equivalent to the infinite suspension spectrum of the classifying
space $\ClSp$, thus partly justifying the use of the term
``classifying spectrum''. In \fullref{sec:Comparison} we extend
this observation to structured classifying spectra.

As an obvious consequence of \fullref{prop:omegaBasis},
the group of homotopy classes of maps between classifying spectra
of fusion systems has an appealingly simple description, analogous
to \fullref{thm:Segal}. This is 
\fullref{mthm:ClSpectraMapBasis} of the introduction.
\begin{theorem} \label{thm:ClSpectraMapBasis}
Let $\F_1$ and $\F_2$ be saturated fusion systems over the finite
$p$--groups $S_1$ and $S_2$, respectively, and let $\wtilde{I}$
be the set of $(\F_1,\F_2)$--conjugacy classes of nontrivial
$(S_1,S_2)$--pairs. Pick a representative $(P_i,\psi_i)$ for each
\mbox{$i \in \wtilde{I}$}. Then the collection
$$\{\sigma_{\F_2} \circ \wtilde{\alpha}([P_i,\psi_i]) \circ t_{\F_1} \mid i\in \wtilde{I} \}  $$
forms a $\Zp$--basis for $\ClSpecMaps{\F_1}{\F_2}$.
\end{theorem}
\begin{proof}
The map
$$\ClSpecMaps{\F_1}{\F_2} \longrightarrow \StableMaps{BS_1}{BS_2},\hspace{.5cm} f\mapsto t_{\F_2} \circ f \circ \sigma_{\F_1}$$
and its left inverse
$$\StableMaps{BS_1}{BS_2} \longrightarrow  \ClSpecMaps{\F_1}{\F_2},\hspace{.5cm} g\mapsto \sigma_{\F_2} \circ g \circ t_{\F_1}$$
make $\ClSpecMaps{\F_1}{\F_2}$ isomorphic to the submodule
$$\StableCharIdem_{\F_2} \circ \StableMaps{BS_1}{BS_2} \circ \StableCharIdem_{\F_1} \subset \StableMaps{BS_1}{BS_2}.$$
Since
$$\StableCharIdem_{\F_2} \circ \StableMaps{BS_1}{BS_2} \circ \StableCharIdem_{\F_1} = \wtilde{\alpha}\left( \CharIdem_{\F_2} \circ \pComp{\wtilde{A}(S_1,S_2)} \circ \CharIdem_{\F_1} \right),$$
the result now follows from the explicit description of the
$\Zp$--basis of \mbox{$\CharIdem_{\F_2} \circ \pComp{A(S_1,S_2)}
\circ \CharIdem_{\F_1}$} given in
\fullref{prop:omegaBasis}.
\end{proof}
This theorem can be applied to fusion systems arising from finite
groups to give a new and simple description of the group of homotopy
classes of stable maps between $p$--completed classifying spaces of
finite groups. Since this result is of independent interest, it will
be presented separately in \cite{KR:NewSegal}. The proof given there
is similar but more direct, with the added advantage that the double
coset formula is preserved, and that the target can be the
classifying space of a compact Lie group.

The following theorem further justifies the use of the term
classifying spectrum. This appears in the introduction as 
\fullref{mthm:sigmaToF}.
\begin{theorem}\label{thm:sigmaToF}
If $\F$ is a saturated fusion system over a finite $p$--group $S$,
then
$$\HomF{P}{Q} = \{\varphi \in \Hom{P}{Q}{}  \mid \sigma_{\F} \circ \ClSpectrum{\iota_Q} \circ \ClSpectrum{\varphi} \simeq \sigma_{\F} \circ \ClSpectrum{\iota_P} \}$$ for all subgroups \mbox{$P,Q \leq S$}.
\end{theorem}
\begin{proof}
When $P$ is the trivial subgroup there is nothing to show, so we
assume otherwise. Now, let \mbox{$\varphi \in \Hom{P}{Q}{}$}.
Since \mbox{$t_{\F} \circ \sigma_{\F} \simeq
\StableCharIdem_{\F}$} and \mbox{$\sigma_{\F} \circ
\StableCharIdem_{\F} \simeq \sigma_{\F}$}, we have
$$\sigma_{\F} \circ \ClSpectrum{\iota_Q} \circ \ClSpectrum{\varphi} \simeq \sigma_{\F} \circ \ClSpectrum{\iota_P} $$
if and only if
\[
  \StableCharIdem_{\F} \circ \ClSpectrum{\iota_Q} \circ \ClSpectrum{\varphi} \simeq \StableCharIdem_{\F} \circ  \ClSpectrum{\iota_P},
\]
which we rewrite as
\[
  \StableCharIdem_{\F} \circ \wtilde{\alpha}([P,\iota_Q \circ \varphi]_P^{S}) \simeq \StableCharIdem_{\F} \circ  \wtilde{\alpha}([P,\iota_P]_P^{S}).
\]
Since $\ClSpectrum{P}$ is naturally equivalent to the classifying
spectrum of $\F_P$, the fusion system of $P$, we can apply 
\fullref{thm:ClSpectraMapBasis} to see that the last equivalence holds
if and only if the $(P,S)$--pairs $(P,\iota_Q \circ \varphi)$ and
$(P,\iota_P)$ are $(\F_P,\F)$--conjugate. By definition this means
that there exist \mbox{$g \in P$} and \mbox{$\varphi' \in
\HomF{P}{\varphi(P)}$} such that the following diagram commutes
\[
\begin{CD}
   {P} @ > \iota_P >> P\\
   @ V \cong V c_g V   @VV \varphi' V \\
   P @>\iota_Q \circ \varphi >> \varphi(P),
\end{CD}
\]
or in other words such that
$$\varphi(x) = \varphi' \circ c_g^{-1}(x)$$
for all \mbox{$x \in P$}. This is in turn true if and only
\mbox{$\varphi \in \HomF{P}{Q}$}.
\end{proof}
This theorem shows in particular that the fusion system of a
finite group $G$ with Sylow subgroup $S$ is determined by
$\ClSpectrum{G}$ regarded as an object under $\ClSpectrum{S}$.
Example 5.2 of \cite{MP3} shows that the homotopy type of
$\ClSpectrum{G}$ alone does not determine the fusion system; the
classifying spectrum must be regarded as an object under
$\ClSpectrum{S}$. As this result is also of independent interest,
and can be proved directly using the Segal Conjecture, it is
treated separately in \cite{KR:NewStableClassification}.
\smallskip

We will show that the assignment of a classifying spectrum to a
saturated fusion systems is functorial, but first we need to
specify which notion of morphisms between fusion systems we are
working with. The following definition appeared in \cite{Puig} but
using different terminology.
\begin{definition} \label{def:FusPresMorphism}
Let $\F_1$ and $\F_2$ be fusion systems over finite $p$--groups
$S_1$ and $S_2$, respectively. A $(\F_1,\F_2)$--fusion-preserving
homomorphism is a group homomorphism \mbox{$\gamma\negmedspace:
S_1 \longrightarrow S_2$} for which there exists a functor
\mbox{$\PresFunctor{\gamma} \negmedspace: \F_1 \longrightarrow
\F_2$} such that
$$\PresFunctor{\gamma}(P) = \gamma(P)$$
for all subgroups \mbox{$P \leq S_1$}, and
$$\gamma\vert_{Q} \circ \varphi = \PresFunctor{\gamma}(\varphi) \circ \gamma\vert_P$$
for all \mbox{$\varphi \in \Hom{P}{Q}{\F_1}$}.
\end{definition}

\begin{remark}
When there is no danger of confusion, we will often say simply that
$\gamma$ is fusion-preserving. The functor $\PresFunctor{\gamma}$ is
uniquely determined by the fusion-preserving morphism $\gamma$ in
the above definition and we will from now on let
$\PresFunctor{\gamma}$ denote the functor defined by a
fusion-preserving homomorphism $\gamma$.
\end{remark}

As a motivation for this definition one may may consider the
equation \mbox{$\gamma \circ c_g(x) = c_{\gamma(g)} \circ
\gamma(x)$} for a group homomorphism \mbox{$\gamma \negmedspace :
G \to H$} and elements \mbox{$g,x \in G$}. Certainly it follows
easily from this equation that a homomorphism between finite
groups induces a fusion-preserving homomorphism between their
fusion systems. Similarly it is not too difficult to prove that a
map between classifying spaces of $p$--local finite groups induces
a fusion-preserving homomorphism between their fusion systems (see
\fullref{rem:UnstableVersion}).

\begin{definition}
Let $\FusSysCat$ be the category whose objects are pairs $(S,\F)$
consisting of a finite $p$--group $S$ and a fusion system $\F$ over
$S$, and whose morphisms are fusion-preserving homomorphisms. Let
$\SatFusSysCat$ be the full subcategory of $\FusSysCat$ whose
objects are pairs $(S,\F)$ where $\F$ is saturated.
\end{definition}

Since a given homomorphism \mbox{$\gamma\negmedspace: S_1 \to
S_2$} can be fusion-preserving for many different fusion systems,
we will write $\gamma_{\F_1}^{\F_2}$ to specify that it is
regarded as an element in the morphism set
$\Mor{(S_1,\F_1)}{(S_2,\F_2)}{\FusSysCat}$.

\begin{lemma}\label{lem:sigmagammaFinv}
Let $\F_1$ and $\F_2$ be saturated fusion systems over finite
$p$--groups $S_1$ and $S_2$, respectively, and let
\mbox{$\gamma\negmedspace : S_1 \to S_2$} be a fusion-preserving
homomorphism. If \mbox{$f\negmedspace : \ClSpectrum{S_2} \to X$}
is a $\F_2$--stable map then
$$f \circ \ClSpectrum{\gamma} \negmedspace : \ClSpectrum{S_1} \to X$$
is $\F_1$--stable, and
$$f \circ \ClSpectrum{\gamma} \circ \StableCharIdem_{\F_1} \simeq f \circ \ClSpectrum{\gamma}.$$
\end{lemma}
\begin{proof}
Indeed, for \mbox{$P \leq S_1$} and \mbox{$\varphi \in
\Hom{P}{S_1}{\F_1}$} we have
\begin{align*}
  f \circ \ClSpectrum\gamma \circ \ClSpectrum{\varphi}
  &\simeq f \circ \ClSpectrum{\left(\gamma \circ \varphi\right)}\\
  &\simeq f \circ \ClSpectrum{\PresFunctor{\gamma}\left(\varphi\right)} \circ \ClSpectrum{\gamma\vert_P}\\
  &\simeq f \circ \ClSpectrum{\iota_{\gamma\left(P\right)}} \circ \ClSpectrum{\gamma}\vert_P\\
  &\simeq f \circ \ClSpectrum{\left(\gamma \circ \iota_P\right)}\\
  &\simeq f \circ \ClSpectrum\gamma \circ \ClSpectrum{\iota_P},
\end{align*}
proving $\F_1$--stability. The second claim follows from 
\fullref{cor:FinvEatsOmega}.
\end{proof}

Now, given saturated fusion systems $\F_1$ and $\F_2$ over finite
$p$--groups $S_1$ and $S_2$, respectively, and a fusion-preserving
homomorphism \mbox{$\gamma\negmedspace: S_1 \to S_2$}, we define a
map of spectra
$$\ClSpectrum{\gamma_{\F_1}^{\F_2}} \negmedspace: \ClSpectrum{\F_1} \to \ClSpectrum{\F_2}$$
by
$$\ClSpectrum{\gamma_{\F_1}^{\F_2}} := \sigma_{\F_2} \circ \ClSpectrum\gamma \circ t_{\F_1}.$$
We will show that this assignment is functorial below. The proof
will first be presented in a setting where the target category keeps
track of the structure maps of the classifying spectra. To this end
we make the following definition.

\begin{definition} \label{def:SCSTargetCat}
Let $\SCSTargetCat$ be the category whose objects are maps of
spectra
$$\sigma\negmedspace: \ClSpectrum{S} \to X,$$
where $S$ is a finite $p$--group and $X$ is a $p$--complete
spectrum, and where the set of morphisms from
\mbox{$\sigma_1\negmedspace: \ClSpectrum{S_1} \to X_1$} to
\mbox{$\sigma_2\negmedspace: \ClSpectrum{S_2} \to X_2$} consists
of pairs $(\gamma,g),$ where \mbox{$\gamma\negmedspace: S_1 \to
S_2$} is a group homomorphism and \mbox{$g\negmedspace: X_1 \to
X_2$} is a map of spectra, such that the following diagram
commutes up to homotopy:
\[
\begin{CD}
{\ClSpectrum{S_1}} @ > {\sigma_1} >> {X_1} \\
@ V {\ClSpectrum{\gamma}} VV @ VV g V \\
{\ClSpectrum{S_2}} @> {\sigma_2} >> {X_2.} \\
\end{CD}
\]
\end{definition}

\begin{theorem} \label{thm:StrClSpecFunctor}
There is a \emph{structured classifying spectrum functor}
$$\StrClSpecFunctor : \SatFusSysCat \longrightarrow \SCSTargetCat$$
defined on objects by
$$(S,\F) \mapsto \left(\sigma_{\F}\negmedspace: \ClSpectrum{S} \to \ClSpectrum{\F}\right)$$
and on morphisms by
$$\gamma_{\F_1}^{\F_2} \mapsto \left(\gamma, \ClSpectrum{\gamma_{\F_1}^{\F_2}}\right).$$
\end{theorem}
\begin{proof}
There are three things to show. First, it is clear by construction
that \mbox{$\sigma_{\F}\negmedspace: \ClSpectrum{S} \to
\ClSpectrum{\F}$} is an object of $\SCSTargetCat$.

Second, if
\mbox{$\gamma_{\F_1}^{\F_2}\in\Mor{(S_1,\F_1)}{(S_2,\F_2)}{\SatFusSysCat}$},
then we  need to prove commutativity of the diagram
\[
\begin{CD}
{\ClSpectrum{S_1}} @ > {\sigma_{\F_1}} >> {\ClSpectrum{\F_1}} \\
@ V {\ClSpectrum\gamma} VV @ VV \ClSpectrum{\gamma_{\F_1}^{\F_2}} V \\
{\ClSpectrum{S_2}} @> {\sigma_{\F_2}} >> {\ClSpectrum{\F_2}.} \\
\end{CD}
\]
Going around the top half of the square, we get
\begin{align*}
  \ClSpectrum{\gamma_{\F_1}^{\F_2}} \circ \sigma_{\F_1}
  &= \sigma_{\F_2} \circ \ClSpectrum\gamma \circ t_{\F_1} \circ \sigma_{\F_1}\\
  &\simeq \sigma_{\F_2} \circ \ClSpectrum\gamma \circ \StableCharIdem_{\F_1}.
\end{align*}
Going around the bottom half of the square, we get
$$\sigma_{\F_2} \circ \ClSpectrum\gamma.$$
Since $\sigma_{\F_2}$ is $\F_2$--stable, these are the same by
\fullref{lem:sigmagammaFinv}.

Third, we need to show that $\StrClSpecFunctor$ preserves
compositions. For this, we let $(S_1,\F_1),(S_2,\F_2)$ and
$(S_3,\F_3)$ be saturated fusion systems, and
$$\gamma_1\negmedspace: S_1 \to S_2 \hspace{1cm} \text{and} \hspace{1cm} \gamma_2\negmedspace: S_2 \to S_3$$
be fusion-preserving morphisms between them. Now,
\begin{align*}
  \ClSpectrum{{\gamma_2}_{\F_2}^{\F_3}} \circ \ClSpectrum{{\gamma_1}_{\F_1}^{\F_2}}
  &\stackrel{\phantom{(\ref{lem:sigmagammaFinv})}}{\simeq}  \sigma_{\F_3} \circ \ClSpectrum{\gamma_2} \circ \underbrace{t_{\F_2} \circ \sigma_{\F_2}}_{=\StableCharIdem_{\F_2}} \circ \ClSpectrum{\gamma_1} \circ t_{\F_1}\\
  &\stackrel{(\ref{lem:sigmagammaFinv})}{\simeq} \sigma_{\F_3} \circ \ClSpectrum{\gamma_2} \circ \ClSpectrum{\gamma_1} \circ t_{\F_1} \\
  &\stackrel{\phantom{(\ref{lem:sigmagammaFinv})}}{\simeq} \sigma_{\F_3} \circ \ClSpectrum{\left(\gamma_2 \circ \gamma_1\right)} \circ t_{\F_1} \\
  &\stackrel{\phantom{(\ref{lem:sigmagammaFinv})}}{\simeq} \ClSpectrum{{\left(\gamma_2 \circ \gamma_1\right)}_{\F_1}^{\F_3}}.
\end{align*}
This completes the proof. \end{proof}

Composing with the forgetful functor from $\SCSTargetCat$ to
$\SpectrumCat$, we obtain the following corollary.
\begin{corollary} \label{cor:ClSpecFunctor}
There is a \emph{classifying spectrum functor}
$$\ClSpectrum : \SatFusSysCat \longrightarrow \SpectrumCat$$
defined on objects by
$$(S,\F) \mapsto \ClSpectrum{\F}$$
and on morphisms by
$$\gamma_{\F_1}^{\F_2} \mapsto \ClSpectrum{\gamma_{\F_1}^{\F_2}}.$$
\end{corollary}

We conclude this section by illustrating that the
fusion-preserving homomorphisms are the only homomorphisms
inducing maps between classifying spectra that preserve structure
maps.
\begin{proposition} Let $\F_1$ and $\F_2$ be saturated fusion
systems over finite \mbox{$p$--groups} $S_1$ and $S_2$, respectively.
If \mbox{$\gamma \negmedspace : S_1 \to S_2$} is a group
homomorphism such that \mbox{$\ClSpectrum\gamma \negmedspace :
\ClSpectrum{S_1} \to \ClSpectrum{S_2}$} restricts to a map \mbox{$g
\negmedspace : \ClSpectrum{\F_1} \to \ClSpectrum{\F_2}$} making the
following diagram commute up to homotopy
\[
\begin{CD}
{\ClSpectrum{S_1}} @ > {\sigma_{\F_1}} >> {\ClSpectrum{\F_1}} \\
@ V {\ClSpectrum\gamma} VV @ VV g V \\
{\ClSpectrum{S_2}} @> {\sigma_{\F_2}} >> {\ClSpectrum{\F_2},} \\
\end{CD}
\]
then $\gamma$ is fusion-preserving.
\end{proposition}
\begin{proof}
We will produce a functor \mbox{$\PresFunctor{\gamma} \negmedspace
: \F_1 \to \F_2 $} that makes $\gamma$ fusion-preserving. There
are two things to check. First, that given a homomorphism
\mbox{$\varphi \in \Hom{P}{Q}{\F_1}$} there is a unique induced
homomorphism \mbox{$\PresFunctor{\gamma}(\varphi) \negmedspace :
\gamma(P) \to \gamma(Q)$} such that \mbox{$\gamma\vert_{Q} \circ
\varphi = \PresFunctor{\gamma}(\varphi) \circ \gamma\vert_P$}, and
second, that $\PresFunctor{\gamma}(\varphi)$ is in $\F_2$.
Functoriality of $\PresFunctor{\gamma}$ follows from the
uniqueness.

To prove the first claim, let $K$ be the kernel of $\gamma$. Then
\mbox{$K \cap P$} is the kernel of $\gamma\vert_P$, and by
standard group theory there exists a homomorphism
$\PresFunctor{\gamma}(\varphi)$ fitting into the following
commutative diagram
\[
\xymatrix{
 {K \cap P} \ar[1,0]  & \\
 {P} \ar[0,1]^{\varphi} \ar@{->>}[1,0]^{\gamma\vert_P} & Q \ar@{->>}[1,0]^{\gamma\vert_Q}\\
 {\gamma(P)} \ar@{-->}[0,1]^{\PresFunctor{\gamma}(\varphi)} &
 {\gamma(Q)}
}
\]
if and only if the restriction of \mbox{$\gamma\vert_Q \circ
\varphi$} to $K \cap P$ is trivial. Furthermore, since
$\gamma\vert_P$ is surjective onto $\gamma(P)$, this condition
uniquely determines $\PresFunctor{\gamma}(\varphi)$ if it exists.
Now,
\begin{align*}
 \sigma_{\F_2} \circ \ClSpectrum{(\gamma\vert_Q \circ \varphi\vert_{K \cap P})}
 &\simeq  \sigma_{\F_2} \circ \ClSpectrum{\gamma} \circ \ClSpectrum{\iota_Q} \circ \ClSpectrum{\varphi} \circ \ClSpectrum{\iota_{K \cap P}}\\
 &\simeq  g \circ \sigma_{\F_1} \circ \ClSpectrum{\iota_Q} \circ \ClSpectrum{\varphi} \circ \ClSpectrum{\iota_{K \cap P}}\\
 &\simeq  g \circ \sigma_{\F_1} \circ \ClSpectrum{\iota_P} \circ \ClSpectrum{\iota_{K \cap P}}\\
 &\simeq  \sigma_{\F_2} \circ \ClSpectrum{\gamma} \circ \ClSpectrum{\iota_{K \cap P}}\\
% &\simeq  \sigma_{\F_2} \circ \ClSpectrum{\gamma\vert_{K \cap P}}\\
 &\simeq  *.
\end{align*}
By applying \fullref{thm:ClSpectraMapBasis} we conclude that
\mbox{$(\gamma\vert_Q \circ \varphi\vert_{K \cap P})$} is trivial.

The second claim is proved similarly by first performing the
following manipulation
\begin{align*}
  \sigma_{\F_2} \circ \ClSpectrum{\iota_{\gamma(Q)}} \circ \ClSpectrum{\PresFunctor{\gamma}(\varphi)} \circ \ClSpectrum{\gamma\vert_P}
 &\simeq \sigma_{\F_2} \circ \ClSpectrum{\iota_{\gamma(Q)}} \circ \ClSpectrum{\gamma\vert_Q} \circ \ClSpectrum{\varphi}\\
 &\simeq  g \circ \sigma_{\F_1} \circ \ClSpectrum{\iota_Q} \circ \ClSpectrum{\varphi}\\
 &\simeq g \circ \sigma_{\F_1} \circ \ClSpectrum{\iota_P}\\
 &\simeq  \sigma_{\F_2} \circ \ClSpectrum{\gamma} \circ \ClSpectrum{\iota_P},
\end{align*}
from which we conclude by \fullref{thm:ClSpectraMapBasis},
that the $(P,S_2)$--pair \mbox{$(P,\iota_{\gamma(Q)} \circ
\PresFunctor{\gamma}(\varphi)\circ {\gamma\vert_P})$} is
$(\F_P,\F_2)$--conjugate to \mbox{$(P,\gamma \circ \iota_P)$}. By
definition this means that there exist \mbox{$g \in P$} and
\mbox{$\varphi' \in \Hom{\gamma(P)}{\gamma(\varphi(P))}{\F_2}$}
making the following diagram commute
\[
\begin{CD}
   {P} @ > \gamma\vert_P >> \gamma(P)\\
   @ V \cong V c_g V   @VV \varphi' V \\
   P @> \PresFunctor{\gamma}(\varphi) \circ \gamma\vert_P >> \gamma(\varphi(P)).
\end{CD}
\]
This implies that
\begin{align*}
  \PresFunctor{\gamma}(\varphi) \circ \gamma (x) &= \varphi' \circ \gamma \circ c_{g^{-1}}(x) \\
  & = \varphi' \circ c_{\gamma(g^{-1})} \circ \gamma(x),
\end{align*}
for all \mbox{$x \in P$}. Since $\gamma$ is surjective onto
$\gamma(P)$ this implies that
$$\PresFunctor{\gamma}(\varphi) = \varphi' \circ c_{\gamma(g^{-1})} \in \Hom{\gamma(P)}{\gamma(Q)}{\F_2}.$$
\end{proof}

\begin{remark} \label{rem:UnstableVersion} \
\begin{enumerate}
  \item
An unstable version of the preceding proof, using
\cite[Proposition 4.4]{BLO2} instead of 
\fullref{thm:ClSpectraMapBasis}, shows that a map between classifying
spaces of $p$--local finite groups restricts to a fusion-preserving
homomorphism of underlying $p$--groups.
  \item In the second paragraph of the proof we showed that
$$\varphi(\operatorname{Ker} \gamma  \cap P) \leq \operatorname{Ker} \gamma$$
for all \mbox{$P,Q \leq S$} and \mbox{$\varphi \in
\Hom{P}{Q}{\F_1}$}. In other words, $\operatorname{Ker} \gamma$ is
\emph{strongly closed in $\F_1$.}
\end{enumerate}
\end{remark}

\section{Transfer theory for classifying spectra} \label{sec:Transfer}
In the classical setting, the map of classifying spaces induced by
an injective homomorphism \mbox{$\psi\negmedspace : G \to H$} of
finite groups admits a transfer \mbox{$tr_{\psi}\negmedspace :
\PtdStable{BH} \to \PtdStable{BG}$} (see for example
\cite{Ad:InfLoopSp}), which restricts to a reduced transfer
\mbox{$\Stable{BH} \to \Stable{BG}$}, also denoted by $tr_\psi$. In
this section we adapt this theory to injective fusion-preserving
homomorphisms and obtain a transfer theory between classifying
spectra. The transfers we construct will not respect structure maps
of saturated fusion systems and so do not fit into the framework of
structured classifying spectra. They do however respect transfer
maps of saturated fusion systems and we will formulate our results
in terms of a functor analogous to the structured classifying space
functor. One motivation for doing so is the following, rather
surprising result.

\begin{proposition}\label{prop:tToF}
If $\F$ is a saturated fusion system over a finite $p$--group $S$,
then
$$\HomF{P}{Q} = \{\varphi \in \Inj{P}{Q} \mid tr_{\varphi} \circ tr_Q \circ t_{\F} \simeq tr_P \circ t_{\F} \}.$$
\end{proposition}
\begin{proof}
The proof is similar to the proof of \fullref{thm:sigmaToF}.
The extra condition that $\varphi$ be injective is required to
ensure the existence of the transfer $tr_\varphi$.
\end{proof}

\begin{definition}
Let $\F_1$ and $\F_2$ be saturated fusion systems over finite
\mbox{$p$--groups} $S_1$ and $S_2$, respectively. If $\gamma$ is an
injective, fusion-preserving homomorphism \mbox{$S_1 \to S_2$}, then
the \emph{$(\F_1,\F_2)$--transfer of $\gamma$} is
$$\Transfer{\gamma_{\F_1}^{\F_2}} := \sigma_{\F_1} \circ tr_{\gamma} \circ t_{\F_2} : \ClSpectrum{\F_2} \to \ClSpectrum{\F_1}.$$
\end{definition}

As in the previous section we show that this assignment is
functorial, working first with a target category that keeps track of
all the structure involved.

\begin{definition} Let $\InjSatFusSysCat$ be the subcategory of
$\SatFusSysCat$, whose objects are saturated fusion systems, but
whose morphisms are the injective fusion-preserving homomorphisms.
\end{definition}

\begin{definition}
Let $\TransferTargetCat$ be the category whose objects are maps of
spectra
$$t\negmedspace: X \to \ClSpectrum{S} ,$$
where $S$ is a finite $p$--group and $X$ is a $p$--complete spectrum,
and where the set of morphisms from \mbox{$t_2\negmedspace:  X_2 \to
\ClSpectrum{S_2}$} to \mbox{$t_1\negmedspace: X_1 \to
\ClSpectrum{S_1}$} consists of pairs $(\gamma,g)$ where
\mbox{$\gamma \negmedspace : S_1 \to S_2$} is a group monomorphism
and \mbox{$g\negmedspace: X_2 \to X_1$} is a map of spectra such
that the following diagram commutes up to homotopy:
\[
\begin{CD}
 {X_2}  @ > {t_2} >> {\ClSpectrum{S_2}}\\
@ V g VV @ VV tr_{\gamma} V \\
{X_1} @> {t_1} >> {\ClSpectrum{S_1}.} \\
\end{CD}
\]
\end{definition}

The proof of functoriality for transfers is very similar to the
proof of \fullref{thm:StrClSpecFunctor} and we will need the
following analogue of \fullref{lem:sigmagammaFinv}.
\begin{lemma}\label{lem:tr_gammatFinv}
Let $\F_1$ and $\F_2$ be saturated fusion systems over finite
$p$--groups $S_1$ and $S_2$, respectively, and let
\mbox{$\gamma\negmedspace : S_1 \to S_2$} be a fusion-preserving
homomorphism. If \mbox{$f\negmedspace : X \to \ClSpectrum{S_2}$}
is a $\F_2$--stable map, then
$$tr_{\gamma} \circ f  \negmedspace : X \to \ClSpectrum{S_1}$$
 is $\F_1$--stable and
$$\StableCharIdem_{\F_1} \circ tr_{\gamma} \circ f  = tr_{\gamma} \circ f. $$
\end{lemma}
\begin{proof}
Let \mbox{$P \leq S_1$} and \mbox{$\varphi \in \Hom{P}{S}{\F_1}$}.
Recalling that the assignment of transfers to group monomorphisms
contravariantly preserves compositions up to homotopy, we get
\begin{align*}
  tr_\varphi \circ tr_\gamma \circ f
  &\simeq tr_{\left(\gamma \circ \varphi\right)} \circ f\\
  &\simeq tr_{\left(\PresFunctor{\gamma}(\varphi)\circ \gamma\vert_P\right)} \circ f\\
  &\simeq tr_{\gamma\vert_P} \circ tr_{\PresFunctor{\gamma}(\varphi)} \circ f  \\
  &\simeq tr_{\gamma\vert_P} \circ tr_{\gamma(P)} \circ f  \\
  &\simeq tr_{\left(\iota_{\gamma(P)} \circ \gamma\vert_P\right)} \circ f \\
  &\simeq tr_{\left(\gamma \circ \iota_P\right)} \circ f\\
  &\simeq tr_P \circ tr_\gamma \circ f,
\end{align*}
proving $\F_1$--stability. The second claim follows from 
\fullref{cor:FinvEatsOmega}.
\end{proof}

\begin{theorem} \label{thm:TransferFunctor}
There is a \emph{structured fusion transfer functor}
$$\TransferFunctor  : \InjSatFusSysCat \to \TransferTargetCat^{op}$$
defined on objects by
$$(S,\F) \mapsto \left(t_{\F}\negmedspace: \ClSpectrum{\F} \to \ClSpectrum{S} \right)$$
and on morphisms by
$$\gamma_{\F_1}^{\F_2} \mapsto \left(\gamma,\Transfer{\gamma_{\F_1}^{\F_2}} \right).$$
\end{theorem}
\begin{proof}
The proof is similar to the proof of 
\fullref{thm:StrClSpecFunctor}. There are again three things to show.
First, it is clear by construction that \mbox{$t_{\F}\negmedspace:
\ClSpectrum{\F} \to \ClSpectrum{S}$} is an object of
$\TransferTargetCat$.

Second, if
\mbox{$\gamma_{\F_1}^{\F_2}\in\Mor{(S_1,\F_1)}{(S_2,\F_2)}{\InjSatFusSysCat}$},
then we prove commutativity of the diagram
\[
\begin{CD}
 {\ClSpectrum{\F_2}}  @ > {t_{\F_2}} >> {\ClSpectrum{S_2}}\\
@ V \Transfer{\gamma_{\F_1}^{\F_2}} VV @ VV tr_{\gamma} V \\
{\ClSpectrum{\F_1}} @> {t_{\F_1}} >> {\ClSpectrum{S_1}.} \\
\end{CD}
\]
by applying \fullref{lem:tr_gammatFinv} and recalling that
$t_{\F_2}$ is $\F_2$--stable (since $\StableCharIdem_{\F_2}$ is
$\F_2$--stable).

The third part, that $\TransferFunctor$ preserves compositions,
also follows from \fullref{lem:tr_gammatFinv} and functoriality
of the classical transfer, just as in the proof of 
\fullref{thm:StrClSpecFunctor}.
\end{proof}

Composing with the forgetful functor from $\TransferTargetCat^{op}$
to $\SpectrumCat^{op}$ we obtain the following corollary.

\begin{corollary}
There is a \emph{fusion transfer functor}
$$\InjSatFusSysCat \to \SpectrumCat^{op}$$
defined on objects by
$$(S,\F) \mapsto \ClSpectrum{\F} $$
and on morphisms by
$$\gamma_{\F_1}^{\F_2} \mapsto \Transfer{\gamma_{\F_1}^{\F_2}}.$$
\end{corollary}

We conclude this section with an application of the transfer
theory. If $\F_1$ and $\F_2$ are saturated fusion systems over a
finite $p$--group $S$, then we say that a map
$$f\negmedspace: \ClSpectrum{\F_2} \to \ClSpectrum{\F_1}$$
is \emph{structure-preserving} if \mbox{$f \circ \sigma_{\F_2}
\simeq \sigma_{\F_1}$}. We say that $\ClSpectrum{\F_2}$ is a
\emph{summand of $\ClSpectrum{\F_1}$ as objects under
$\ClSpectrum{S}$} if there is a structure-preserving map
\mbox{$f\negmedspace: \ClSpectrum{\F_1} \to \ClSpectrum{\F_2}$}
and a map \mbox{$r\negmedspace: \ClSpectrum{\F_2} \to
\ClSpectrum{\F_1} $} such that
$$f \circ r \simeq id_{\ClSpectrum{\F_2}}.$$

\begin{proposition}
Let $\F_1$ and $\F_2$ be saturated fusion systems over a finite
\mbox{$p$--group $S$}. The following four conditions are
equivalent:
\begin{enumerate}
  \item[(i)] $\F_1$ is a subcategory of $\F_2$;
  \item[(ii)] $id_S$ is $(\F_1,\F_2)$--fusion-preserving;
  \item[(iii)] $\ClSpectrum{\F_2}$ is a summand of
  $\ClSpectrum{\F_1}$ as objects under $\ClSpectrum{S}$;
  \item[(iv)] There is a structure-preserving map \mbox{$f\negmedspace: \ClSpectrum{\F_1}\to
  \ClSpectrum{\F_2}$}.
\end{enumerate}
\end{proposition}
\begin{proof} (i) $\Rightarrow$ (ii): Obvious.\\
(ii) $\Rightarrow$ (iii): The map
\mbox{$\ClSpectrum{id_{\F_1}^{\F_2}}$} is structure-preserving and
\begin{align*}
  \ClSpectrum{id_{\F_1}^{\F_2}} \circ \Transfer{id_{\F_1}^{\F_2}}
  &\stackrel{\phantom{((\ref{lem:sigmagammaFinv})}}{=} \sigma_{\F_2} \circ \ClSpectrum{ id} \circ t_{\F_1} \circ \sigma_{\F_1} \circ tr_{id} \circ t_{\F_2}\\
  &\stackrel{\phantom{(\ref{lem:sigmagammaFinv})}}{=} \sigma_{\F_2} \circ id \circ \StableCharIdem_{\F_1} \circ id \circ t_{\F_2}\\
  &\stackrel{(\ref{lem:sigmagammaFinv})}{=} \sigma_{\F_2} \circ t_{\F_2}\\
  &\stackrel{\phantom{(\ref{lem:sigmagammaFinv})}}{=} id_{\ClSpectrum{\F_2}}.
\end{align*}
(iii) $\Rightarrow$ (iv): Obvious.\\
(iv) $\Rightarrow$ (i): Let \mbox{$P \leq S$} and \mbox{$\varphi
\in \Hom{P}{S}{\F_1}$}. Then
$$\sigma_{\F_1} \circ \ClSpectrum\varphi \simeq \sigma_{\F_1} \circ \ClSpectrum\iota_P.$$
Consequently
$$\sigma_{\F_2} \circ \ClSpectrum\varphi \simeq f \circ \sigma_{\F_1} \circ \ClSpectrum\varphi \simeq f \circ \sigma_{\F_1} \circ \ClSpectrum\iota_P \simeq \sigma_{\F_2} \circ \ClSpectrum\iota_P.$$
By \fullref{thm:sigmaToF} it follows
that \mbox{$\varphi \in \Hom{P}{S}{\F_2}$}.
\end{proof}

\section{Cohomology of classifying spectra} \label{sec:Cohomology}
In this section we apply the \mbox{$mod~ p$} cohomology functor to
the theory of the previous sections and observe that the
cohomology of fusion systems and their classifying spectra behaves
like group cohomology in many important ways. Cohomology will
always be taken with $\Fp$--coefficients and therefore we will
denote the functor $\Coh{-;\Fp}$ by $\Coh{-}$.

We begin by recording the analogue of the last part of 
\fullref{prop:LWImportance}
\begin{proposition} \label{prop:AlsoHolds}
Let $\F$ be a saturated fusion system over a finite $p$--group $S$.
Then the map $\StableCharIdem_{\F}^*$ induced by its stable
characteristic idempotent in cohomology is an idempotent in
$End(\Coh{BS})$, is $\Coh{\F}$--linear and a homomorphism of modules
over the Steenrod algebra; and
$$Im[\Coh{BS}\stackrel{\StableCharIdem_{\F}^*}{\longrightarrow}\Coh{BS}] = \Coh{\F}.$$
\end{proposition}
\begin{proof}
The proof is essentially the same as in \cite{BLO2}. The reader is
referred to \cite[Proposition 5.5]{BLO2} for details.
\end{proof}

The next proposition, which is arguably the most important result
in this section, has already been observed by Broto--Levi--Oliver in
\cite{BLO2}.
\begin{proposition}{\rm\cite[Section 5]{BLO2}}\qua
Let $\F$ be a saturated fusion system over a finite $p$--group $S$.
Then the structure map $\sigma_{\F}$ induces a split monomorphism
$$\Coh{\ClSpectrum\F} \hookrightarrow \Coh{BS}$$
in cohomology, with image $\Coh{\F}$.
\end{proposition}
\begin{proof}
By construction it is clear that $\sigma_{\F}$ induces the
inclusion
$$Im[\Coh{BS}\stackrel{\StableCharIdem_{\F}^*}{\longrightarrow}\Coh{BS}] \hookrightarrow \Coh{BS}$$
in cohomology, and by \fullref{prop:AlsoHolds} this has
image $\Coh{\F}$.
\end{proof}

\begin{remark} Using \fullref{cor:FinvEatsOmega},
Castellana--Morales show in \cite{CM} that the analogous result holds
for any generalized cohomology theory. That is, for a saturated
fusion system $\F$ over a finite $p$--group $S$ and a generalized
cohomology theory $E$, the structure map $\sigma_{\F}$ induces a
split monomorphism
$$E^*(\ClSpectrum\F) \hookrightarrow E^*(BS)$$
$$E^*(\F) := \invlimlim{\F}E^*(B(-)).\leqno{\hbox{with image}}
$$
\end{remark}

Since $\Coh{\F}$ is a subring of $\Coh{BS}$, and is consequently an
unstable algebra over the Steenrod algebra, we can now regard
$\Coh{\ClSpectrum{\F}}$ as an unstable algebra over the Steenrod
algebra via the isomorphism
$$\sigma_{\F}^* : \Coh{\ClSpectrum{\F}} \stackrel{\cong}{\longrightarrow} \Coh{\F}.$$
With this understanding, the map
$$\sigma_{\F}^* : \Coh{\ClSpectrum{\F}} {\longrightarrow} \Coh{\ClSpectrum{S}}$$
is obviously a map of unstable algebras over the Steenrod algebra.

The $\Coh{\F}$--linearity of \fullref{prop:AlsoHolds} is a
version of what is commonly referred to as Frobenius reciprocity
in the classical group cohomology setting. We record this property
in the more familiar form.
\begin{proposition}[Frobenius reciprocity]
Let $\F$ be a saturated fusion system over a finite $p$--group $S$.
If \mbox{$x \in \Coh{\ClSpectrum{\F}}$} and \mbox{$y \in \Coh{BS}$},
then
$$t_{\F}^*\left(\sigma_{\F}^*\left(x\right) \cdot y\right) = x \cdot t_{\F}^*\left(y\right).$$
\end{proposition}
\begin{proof}
Since there is only one fusion system in play, we will drop the
subscript $\F$ and simply write $\sigma,$ $t$ and
$\StableCharIdem$. Now, by using $\Coh{\F}$--linearity of
$\StableCharIdem^*$ and the fact that $\sigma^*$ preserves
multiplication, we have
\begin{align*}
  t^* \left(\sigma^*\left(x\right) \cdot y\right) &= t^* \circ \StableCharIdem^* (\underbrace{\sigma^*\left(x\right)}_{\in \Coh{\F}} \cdot~ y)\\
                            &= t^* \left(\sigma^* \left(x \right) \cdot \StableCharIdem^*\left(y\right)\right)\\
                            &= t^* \left(\vphantom{\sum_{i}} \sigma^*\left(x\right) \cdot \left(\sigma^*  \circ t^*\left(y\right)\right)\right)\\
                            &= \underbrace{t^* \circ \sigma^*}_{id} \left(x \cdot t^*\left(y\right)\right)\\
                            &= x \cdot t^*\left(y\right).
\end{align*}
This completes the proof. \end{proof}

The usefulness of transfers in group cohomology lies to a large
extent in the well known result, that for a group monomorphism
\mbox{$\gamma\negmedspace: G \to H$}, the  effect of the
composition \mbox{$tr_{\gamma}^* \circ \ClSpectrum\gamma^*$} in
cohomology is simply multiplication by the index $[H:\gamma(G)]$.
The situation for transfers between classifying spectra induced by
fusion-preserving monomorphisms is similar. Although one might
expect that the order of the outer automorphism groups of the
fusion systems should come into play, it is only the order of the
underlying $p$--groups that is important. The reason is that the
characteristic idempotents, which are used to construct the
transfer, have a ``normalizing effect''. In some sense they divide
out the order of the outer automorphism groups.
\begin{proposition} \label{prop:ActsByMultiplication}
Let $\F_1$ and $\F_2$ be saturated fusion systems over $S_1$ and
$S_2$, respectively, and let \mbox{$\gamma\negmedspace: S_1 \to
S_2$} be a fusion-preserving monomorphism. Then the composition
$$\Transfer{\gamma_{\F_1}^{\F_2}}^* \circ \ClSpectrum{\gamma_{\F_1}^{\F_2}}^*$$
acts on $\Coh{\ClSpectrum{\F_1}}$ as multiplication by
$|S_2|/|S_1|$.
\end{proposition}
\begin{proof}
Indeed,
\begin{align*}
  \ClSpectrum{\gamma_{\F_1}^{\F_2}} \circ \Transfer{\gamma_{\F_1}^{\F_2}}
  &\stackrel{\phantom{(\ref{lem:sigmagammaFinv})}}{=} \sigma_{\F_2} \circ B\gamma \circ t_{\F_1} \circ \sigma_{\F_1} \circ tr_{\gamma} \circ t_{\F_2} \\
  &\stackrel{\phantom{(\ref{lem:sigmagammaFinv})}}{=} \sigma_{\F_2} \circ B\gamma \circ \StableCharIdem_{\F_1} \circ tr_{\gamma} \circ t_{\F_2} \\
  &\stackrel{(\ref{lem:sigmagammaFinv})}{=} \sigma_{\F_2} \circ B\gamma \circ tr_{\gamma} \circ t_{\F_2}.
\end{align*}
Now, if \mbox{$x \in \Coh{\ClSpectrum{\F_1}},$} then
\begin{align*}
  \Transfer{\gamma_{\F_1}^{\F_2}}^* \circ \ClSpectrum{\gamma_{\F_1}^{\F_2}}^*\left(x\right)
  &= t_{\F_2}^* \circ tr_{\gamma}^* \circ B\gamma^* \circ \sigma_{\F_2}^* \left(x\right) \\
  &= t_{\F_2}^* \left(\frac{|S_2|}{|S_1|} \cdot   \sigma_{\F_2}^* \left(x\right)\right) \\
  &= \frac{|S_2|}{|S_1|} \cdot   t_{\F_2}^* \left(\sigma_{\F_2}^* \left(x\right)\right) \\
  &= \frac{|S_2|}{|S_1|} \cdot x.
\end{align*}
This completes the proof. \end{proof}

Transfers between classifying spectra also exhibit Frobenius
reciprocity.
\begin{proposition}[Frobenius reciprocity]\label{prop:TransfersFrobenius}
Let $\F_1$ and $\F_2$ be saturated fusion systems over $S_1$ and
$S_2$, respectively, and let \mbox{$\gamma\negmedspace: S_1 \to
S_2$} be a fusion-preserving monomorphism. If \mbox{$x \in
\Coh{\ClSpectrum{\F_2}}$} and \mbox{$y \in
\Coh{\ClSpectrum{\F_1}}$}, then
$$\Transfer{\gamma_{\F_1}^{\F_2}}^*\left(\ClSpectrum{\gamma_{\F_1}^{\F_2}}^*\left(x\right) \cdot y\right) = x \cdot \Transfer{\gamma_{\F_1}^{\F_2}}^*\left(y\right).$$
\end{proposition}
\begin{proof}
Using the facts that $\sigma_{\F_1}^*$ is a ring homomorphism, and
that $t_{\F_2}^*$ and $tr_{\gamma}^*$ both exhibit Frobenius
reciprocity, we get
\begin{align*}
  \Transfer{\gamma_{\F_1}^{\F_2}}^* \left(\ClSpectrum{\gamma_{\F_1}^{\F_2}}^*\left(x\right) \cdot y \vphantom{\sum_i}\right)
  &= t_{\F_2}^* \circ tr_{\gamma}^* \circ \sigma_{\F_1}^* \left(\ClSpectrum{\gamma_{\F_1}^{\F_2}}^*\left(x\right) \cdot y \vphantom{\sum_i}\right)\\
  &= t_{\F_2}^* \circ tr_{\gamma}^* \circ \left((\sigma_{\F_1}^* \circ \ClSpectrum{\gamma_{\F_1}^{\F_2}}^*\left(x\right)) \cdot \sigma_{\F_1}^*\left(y\right) \vphantom{\sum_i}\right)\\
  &= t_{\F_2}^* \circ tr_{\gamma}^* \left((\ClSpectrum{\gamma}^* \circ \sigma_{\F_2}^*\left(x\right)) \cdot \sigma_{\F_1}^*\left(y\right) \vphantom{\sum_i}\right)\\
  &= t_{\F_2}^* \left(\sigma_{\F_2}^*\left(x\right) \cdot (tr_{\gamma}^* \circ \sigma_{\F_1}^*\left(y\right)) \vphantom{\sum_i}\right)\\
  &= x \cdot (t_{\F_2}^* \circ tr_{\gamma}^* \circ \sigma_{\F_1}^*\left(y\right)) \vphantom{\sum_i}\\
  &= x \cdot \Transfer{\gamma_{\F_1}^{\F_2}}^*\left(y\right).
\end{align*}
This completes the proof. \end{proof}

The reciprocity results in this section can be shown to hold at
the level of stable homotopy, as is the case with transfers
induced by finite covers, when the structure maps preserve
diagonals (in particular the classifying spectra involved must
have diagonal maps). Since the only known cases where this happens
is for fusion systems belonging to $p$--local finite groups, in
which case the classifying spectra are suspension spectra, this
discussion is postponed for \cite{KR:Transfers-plfgs}.

\section{Comparison to stable classifying spaces} \label{sec:Comparison}
We conclude the paper by comparing the theory of classifying
spectra of saturated fusion systems with the theory obtained by
infinite suspension of classifying spaces of finite groups and
$p$--local finite groups, and proving that the theory of
classifying spectra extends both these theories.

It is shown in \cite{BLO2} that when $\F$ has an associated centric
linking system $\Link$, the classifying spectrum $\ClSpectrum{\F}$
is homotopy equivalent to the $p$--completed suspension spectrum
$\pComp{\Stable{|\Link|}}$. We extend this observation to structured
classifying spectra. The reader is referred to \cite{BLO2} for the
precise definition of centric linking systems and $p$--local finite
groups.
\begin{proposition} \label{prop:SuspendPLFGs}
Let $(S,\F,\Link)$ be a $p$--local finite group. Then the infinite
suspension
$$\Stable{\theta}\negmedspace: \Stable{BS} \longrightarrow \Stable{\ClSp}$$
of the natural inclusion \mbox{$\theta \negmedspace : BS  \to
\ClSp$} is equivalent to the structure map
$$\sigma_{\F}\negmedspace : \ClSpectrum{S} \longrightarrow \ClSpectrum{\F}.$$
\end{proposition}
\begin{proof}
First we recall that $\ClSpectrum{S}$ is homotopy equivalent to
$\Stable{BS}$, and we may therefore identify the two via a chosen
homotopy equivalence.

In \cite{BLO2} it is shown that $\theta$ is $\F$--stable and it
follows that $\Stable{\theta}$ is $\F$--stable. By 
\fullref{cor:FinvEatsOmega} it follows that
$$\Stable{\theta} \circ \StableCharIdem_\F \simeq \Stable{\theta}.$$
By construction of the structured classifying spectrum as a
mapping telescope we get a map \mbox{$h\negmedspace :
\ClSpectrum{\F} \to \Stable{\ClSp}$} such that
$$h \circ \sigma_{\F} \simeq \Stable{\theta}.$$
In particular the corresponding equality holds for the induced maps
in cohomology with $\Fp$--coefficients. It is shown in \cite{BLO2}
(see also \fullref{prop:AlsoHolds}) that in cohomology the
maps $\sigma_{\F}$ and $\theta$ both induce injctions with image
$\Coh{\F}$ in $\Coh{BS}$, and therefore we conclude that $h$ induces
an isomorphism
$$h^* \negmedspace : \Coh{\Stable{\ClSp}} \stackrel{\cong}{\longrightarrow} \Coh{\ClSpectrum{\F}}.$$
Since the spectra involved are $p$--complete, we deduce that $h$
is a homotopy equivalence.
\end{proof}

We now turn our attention to $p$--completed classifying spaces of
finite groups. This theory overlaps in parts with the theory of
classifying spaces of $p$--local finite groups since the
classifying space of the $p$--local finite group induced by a
finite group $G$ is homotopy equivalent to $\pComp{BG}$. An
additional aspect for the stable classifying spaces of groups is
that the inclusion of a Sylow subgroup $S$ into a finite group $G$
has a stable transfer, which we compare with the transfer of a
saturated fusion system.
\begin{proposition} \label{prop:SuspendGroups}
Let $G$ be a finite group with Sylow subgroup $S$. Then the map
$$\ClSpectrum{\iota_{S}}\negmedspace : \ClSpectrum{S} \longrightarrow \ClSpectrum{G}$$
induced by the inclusion \mbox{$S \leq G$} is equivalent to the
structure map
$$\sigma_{\F_S(G)} \negmedspace: \ClSpectrum{S} \longrightarrow \ClSpectrum{\F_S(G)}.$$
Furthermore, if we let $g$ be a homotopy inverse of the homotopy
equivalence
$$\ClSpectrum{\iota_S} \circ tr_S \negmedspace : \ClSpectrum{G} \longrightarrow \ClSpectrum{G},$$
then the map
$$t' := tr_S \circ g \negmedspace : \ClSpectrum{G} \longrightarrow \ClSpectrum{S} $$
is equivalent to the map
$$t_{\F_S(G)} \negmedspace: \ClSpectrum{\F_S(G)} \longrightarrow \ClSpectrum{S}.$$
\end{proposition}
\begin{proof}
Write \mbox{$\F := \F_S(G)$}. It was shown in \cite{BLO2} that for a
finite $p$--group $G$ with Sylow subgroup $S$, the map
\mbox{$B\iota_S \negmedspace : BS \longrightarrow BG $} is
equivalent to the inclusion \mbox{$\theta \negmedspace : BS
\longrightarrow \ClSp$} of $BS$ into the classifying space of the
corresponding $p$--local finite group. By 
\fullref{prop:SuspendPLFGs} it follows that the map
\mbox{$\ClSpectrum{\iota_{S}}\negmedspace : \ClSpectrum{S}
\longrightarrow \ClSpectrum{G}$} is equivalent to the map
\mbox{$\sigma_{\F} \negmedspace: \ClSpectrum{S} \longrightarrow
\ClSpectrum{\F}$}.

When proving the second claim we can fix a homotopy equivalence
\mbox{$\ClSpectrum{\F} \to \ClSpectrum{G}$} and regard $t'$ as a
map $\ClSpectrum{\F} \to \ClSpectrum{S}$. Assume, for now, that
\mbox{$\StableCharIdem_\F \circ tr_S \simeq tr_S.$} Then
\mbox{$\StableCharIdem_\F \circ t' \simeq t'$} and we get
$$t_\F \simeq t_\F \circ \ClSpectrum{\iota_S} \circ t' \simeq t_\F \circ \sigma_\F \circ t' \simeq \StableCharIdem_\F \circ t' \simeq t',$$
which is what we want to show.

To prove that \mbox{$\CharIdem_\F \circ tr_S \simeq tr_S$} it
suffices, by \fullref{cor:FinvEatsOmega}, to establish that
$tr_S$ is $\F$--stable, which is actually quite well known. One way
to convince oneself of this is to note that $tr_S$ is the image in
$\StableMaps{BG}{BS}$ of \mbox{$[S,id_S] \in A(G,S)$}, and that for
\mbox{$P \leq S$} and \mbox{$\varphi \in \HomF{P}{S}$} one has
$$[\varphi(P),\varphi^{-1}]_S^P \circ [S,id_S]_G^S = [\varphi(P),\varphi^{-1}]_G^P = [P,id_P]_G^P \in A(G,P),$$
since $\varphi$ is a conjugation induced by an element of $G$.
\end{proof}

With the notation of the preceding theorem, the stable
characteristic idempotent of $\F_S(G)$ can be obtained as
$$\StableCharIdem_{\F_S(G)} \simeq t' \circ \ClSpectrum{\iota_S},$$
and the characteristic idempotent is then
$$\CharIdem_{\F_S(G)} = \wtilde{\alpha}^{-1}(\StableCharIdem_{\F_S(G)}),$$
regarded as an element of $\pComp{A(S,S)}$.

Another feature of the theory of $p$--completed classifying spaces of
finite groups is that it is functorial. Namely, given a homomorphism
\mbox{$\bar{\gamma} \negmedspace : G_1 \to G_2$} between finite
groups, one gets a map of $p$--completed classifying spaces
\mbox{$\pComp{B\bar{\gamma}} \negmedspace : \pComp{BG_1} \to
\pComp{BG_2}$}. Furthermore, the restriction \mbox{$\gamma
\negmedspace : S_1 \to S_2$} to Sylow subgroups is fusion preserving
for the fusion systems $\F_{S_1}(G_1)$ and $\F_{S_2}(G_2)$. Hence we
get a map \ClSpectrum{\gamma_{\F_1}^{\F_2}} of classifying spectra,
which we can compare to $\ClSpectrum{\bar{\gamma}}$, the infinite
suspension of $\pComp{B\bar{\gamma}}$.
\begin{proposition}
Let \mbox{$\bar{\gamma} \negmedspace : G_1 \to G_2$} be a
homomorphism of finite groups with a restriction \mbox{$\gamma
\negmedspace : S_1 \to S_2$} to Sylow subgroups. Then the diagram
\[
\begin{CD}
{\ClSpectrum{S_1}} @ > {\ClSpectrum{\iota_{S_1}}} >> {\ClSpectrum{G_1}} \\
@ V {\ClSpectrum{\gamma}} VV @ VV \ClSpectrum{\bar\gamma} V \\
{\ClSpectrum{S_2}} @> {\ClSpectrum{\iota_{S_2}}} >> {\ClSpectrum{G_2}} \\
\end{CD}
\]
is equivalent to the diagram
\[
\begin{CD}
{\ClSpectrum{S_1}} @ > \sigma_{\F_1} >> {\ClSpectrum{\F_{S_1}(G_1)}} \\
@ V {\ClSpectrum{\gamma}} VV @ VV \ClSpectrum{\gamma_{\F_1}^{\F_2}} V \\
{\ClSpectrum{S_2}} @> \sigma_{\F_2} >> {\ClSpectrum{\F_{S_2}(G_2)}.} \\
\end{CD}
\]
\end{proposition}
\begin{proof}
By \fullref{prop:SuspendGroups} we can replace all but the
map on the right side of the upper diagram with the corresponding
maps and objects on the lower diagram. All that remains is to show
is now that, when thus regarded as a map
\mbox{${\ClSpectrum{\F_{S_1}(G_1)}} \to
{\ClSpectrum{\F_{S_2}(G_2)}}$}, the map $\ClSpectrum{\bar\gamma}$
is homotopic to the map $\ClSpectrum{\gamma_{\F_1}^{\F_2}}$. Now,
from the equivalence \mbox{$\ClSpectrum{\bar\gamma} \circ
\sigma_{\F_1} \simeq \sigma_{\F_2} \circ \ClSpectrum{\gamma}$} it
follows that
$$\ClSpectrum{\bar\gamma} \simeq \ClSpectrum{\bar\gamma} \circ \sigma_{\F_1} \circ t_{\F_1} \simeq \sigma_{\F_2} \circ \ClSpectrum{\gamma} \circ t_{\F_1} = \ClSpectrum{\gamma_{\F_1}^{\F_2}}.\proved$$
\end{proof}

\bibliographystyle{gtart} \bibliography{link}

\end{document}